\documentclass{article}
\usepackage{amssymb, epsfig}

\newtheorem{theorem}{Theorem}
\newtheorem{lemma}{Lemma}
\newtheorem{proposition}{Proposition}
\newtheorem{corollary}{Corollary}
\newtheorem{definition}{Definition}
\newtheorem{condition}{Condition}
\newtheorem{remark}{Remark}
\newtheorem{example}{Example}

\newcommand{\na}{n}
\newcommand{\nl}{{N-n}}
\newcommand{\bnl}{{(N-n)}}

\newcommand{\CC}{{\mathbb C }}
\newcommand{\EE}{{\mathbb E }}

\newcommand{\LL}{{\mathbb L }}
\newcommand{\MM}{{\mathbb M }}

\newcommand{\PP}{{\mathbb P }}
\newcommand{\QQ}{{\mathbb Q }}
\newcommand{\RR}{{\mathbb R }}
\newcommand{\TT}{{\mathbb T }}
\newcommand{\ZZ}{{\mathbb Z }}
\newcommand{\TTT}{\widetilde{\mathbb T }}

\newcommand{\cA}{{\mathcal A }}
\newcommand{\cB}{{\mathcal B }}
\newcommand{\cC}{{\mathcal C }}
\newcommand{\cD}{{\mathcal D }}
\newcommand{\cE}{{\mathcal E }}
\newcommand{\cF}{{\mathcal F }}

\newcommand{\cJ}{{\mathcal J }}

\newcommand{\cL}{{\mathcal L }}
\newcommand{\cM}{{\mathcal M }}
\newcommand{\cO}{{\mathcal O }}
\newcommand{\cP}{{\mathcal P }}
\newcommand{\cR}{{\mathcal R }}
\newcommand{\cS}{{\mathcal S }}
\newcommand{\cU}{{\mathcal U }}
\newcommand{\cV}{{\mathcal V }}
\newcommand{\cW}{{\mathcal W }}

\newcommand{\gT}{{\mathcal T }\!}

\newcommand{\gB}{{\mathfrak B }}
\newcommand{\gD}{{\mathfrak D }}
\newcommand{\gP}{{\mathfrak P }}
\newcommand{\ga}{{\mathfrak a }}
\newcommand{\gb}{{\mathfrak b }}
\newcommand{\gu}{{\mathfrak u }}

\newcommand{\sA}{{\mathsf A }}
\newcommand{\sB}{{\mathsf B }}
\newcommand{\sZ}{{\mathsf Z }}

\newcommand{\ec}{{\mathsf c }}
\newcommand{\ed}{{\mathsf d }}
\newcommand{\ee}{{\mathsf e }}
\newcommand{\ek}{{\mathsf k }}
\newcommand{\ep}{{\mathsf p }}
\newcommand{\es}{{\mathsf s }}

\newcommand{\eu}{{\mathsf u }}
\newcommand{\ev}{{\mathsf v }}
\newcommand{\ew}{{\mathsf w }}
\newcommand{\ex}{{\mathsf x }}
\newcommand{\ey}{{\mathsf y }}
\newcommand{\ez}{{\mathsf z }}

\newcommand{\sfm}{{\mathsf m }}

\newcommand{\id}{{\mathsf 1 }}

\newcommand{\tsZ}{\widetilde{\mathsf Z }}
\newcommand{\tOmega}{\widetilde{\Omega}}

\newcommand{\conv}{\mathrm{ conv}\,}
\newcommand{\core}{\mathrm{ core}\,}
\newcommand{\Ann}{\mathrm{ Ann}\,}
\newcommand{\rank}{\mathrm{ rank}\,}
\newcommand{\spec}{\mathrm{ spec}\,}
\newcommand{\vol}{\mathrm{ vol}\,}
\newcommand{\rel}{\,\mathrm{ rel }\,}

\newcommand{\tp}{{2\pi{\mathsf i}\, }}

\newcommand{\qed}{{\hfill$\boxtimes$}}

\newcommand{\Hom}[2]{\mathrm{ Hom}(#1,\,#2)}

\newcommand{\modquot}[2]{\mbox{\raisebox{.1em}{$#1$}\hspace{-2mm}
{ / }\hspace{-2mm} \raisebox{-.1em}{$#2$}}}

\begin{document}

\title{Resonant Hypergeometric Systems\\ and Mirror Symmetry }

\author{ Jan Stienstra \\
\footnotesize{Mathematisch Instituut, Universiteit Utrecht}\\
\footnotesize{Postbus 80.010,3508 TA Utrecht, The Netherlands}\\
\footnotesize{E-mail: stien@math.ruu.nl}}

\date{}

\maketitle
\begin{abstract}
In Part I the $\Gamma$-series of \cite{gkz1} are adapted
so that they give solutions for certain resonant systems of
Gel'fand-Kapranov-Zelevinsky hypergeometric differential equations.
For this some complex parameters in the $\Gamma$-series are replaced by
nilpotent elements from a ring
$\cR_{\sA,\gT}\,.$ The adapted $\Gamma$-series is a function
$\Psi_{\gT,\beta}$ with values in the finite dimensional vector space
$\cR_{\sA,\gT}\otimes_\ZZ\CC\,.$
Part II describes applications of these results in the context of toric
Mirror Symmetry. Building on Batyrev's work  \cite{bat2} we show that
a certain relative cohomology module
$H^\na(\TTT\rel\tsZ_{\es-1})$ is a GKZ hypergeometric $\cD$-module
which over an appropriate domain is isomorphic to the trivial
$\cD$-module $\cR_{\sA,\gT}\otimes\cO_\gT$, where $\cO_\gT$ is the sheaf of
holomorphic functions on this domain. The isomorphism is explicitly given by
adapted $\Gamma$-series. As a result one finds the periods of
a holomorphic differential form of degree $d$ on a $d$-dimensional
Calabi-Yau manifold, which are needed for the B-model side input to Mirror
Symmetry. Relating our work with that of Batyrev and Borisov  \cite{babo}
we interpret the ring $\cR_{\sA,\gT}$ as the cohomology ring of a toric variety
and a certain principal ideal in it as a subring of
the Chow ring of a Calabi-Yau complete intersection.
This interpretation takes place on the A-model side of
Mirror Symmetry.
\end{abstract}


\begin{tabbing}
\textbf{Contents}\= \hspace{6mm}\=  \\
Part I \> \> Introduction I \\
\> \ref{regtrifan}.\>
Regular triangulations and the pointed secondary fan \\
\> \ref{thering}.\>
The ring $\cR_{\sA,\gT}\,.$ \\
\> \ref{solutions}.\>
A domain of definition for the function $\Psi_{\gT,\beta}\,.$ \\
\> \ref{betanul}.\>
The special case $\beta=0$\\
\> \ref{coresection}.\>
Triangulations with non-empty core\\
Part II \> \>Introduction II \\
Part IIB\> \>Introduction IIB \\
\> \ref{vmhs}.\>
VMHS associated with a Gorenstein cone \\
Part IIA\> \>Introduction IIA\\
\> \ref{gorenstein}.\>
Triangulations with non-empty core and \\
\> \> completely split reflexive Gorenstein cones \\
\> \ref{tritorvar}. \>
Triangulations and toric varieties \\
\> \ref{cicy}.\>
Calabi-Yau complete intersections in toric varieties\\
Conclusions\\
References
\end{tabbing}

\

\

\begin{center}\textbf{PART I}\end{center}

\section*{Introduction I}

A GKZ hypergeometric system \cite{gkz1} depends on four parameters:
two positive integers $N$ and $\na$, a set
$\{\ga_1,\ldots,\ga_N\}$ of vectors
in $\ZZ^\na$ and a vector $\beta$ in $\CC^{\na}$.
The standard assumptions \cite{gkz1} are

\

\begin{condition}\label{gkzcond}
\footnote{ ${}\;\;\;\ZZ^\na,\RR^\na,\CC^\na$ resp.
$\ZZ^{\na\vee},\RR^{\na\vee},\CC^{\na\vee}$
denote spaces of column vectors resp. row vectors.}\vspace{-7mm}
\begin{eqnarray}
\label{eq:aspan}
&& \ga_1,\ldots,\ga_N \textrm{ generate a rank }n \textrm{ sub-lattice } \MM
\textrm{ in }\ZZ^\na \\
\label{eq:hyperplane}
&& \exists\;\ga_0^\vee\in\ZZ^{\na\vee} \textrm{ such that }
\ga_0^\vee\cdot\ga_i \,=\, 1 \;\;(i=1,\ldots,N)
\end{eqnarray}
\end{condition}
The GKZ system with these parameters is the following system of partial
differential equations for functions $\Phi$ on a torus with coordinates
$v_1,\ldots,v_N$:
\begin{eqnarray}
\label{eq:gkzlin}
\left(-\beta+
\sum_{j=1}^N\:\ga_j\,v_j\frac{\partial}{\partial v_j}\right)\;\Phi &=&0\\
\label{eq:gkzmon}
\left(\prod_{\ell_j>0} \left[\frac{\partial}{\partial v_j}\right]^{\ell_j}
\:-\:\prod_{\ell_j<0} \left[\frac{\partial}{\partial v_j}\right]^{-\ell_j}
\right)\Phi&=&0\hspace{3mm}\textrm{ for } \ell\in\LL
\end{eqnarray}
where (\ref{eq:gkzlin}) is in fact a system of $\na$ equations and
\begin{equation}\label{eq:L}
\LL\::=\:\{\;\ell=(\ell_1,\ldots,\ell_N)^t\in\ZZ^N\:|\:
\ell_1\ga_1+\ldots+\ell_N\ga_N\,=\,0\;\}\,.
\end{equation}
Some of the above data are displayed in the following short exact sequence
in which $\cA$ denotes the linear map $\cA\,:\;\ZZ^N\rightarrow\ZZ^\na$,
$\cA(\lambda)\:=\:\lambda_1\ga_1+\ldots+\lambda_N\ga_N$.
\begin{equation}\label{eq:ses}
0\rightarrow\LL\longrightarrow\ZZ^N\stackrel{\cA}{\longrightarrow}
\MM\rightarrow 0
\end{equation}
\emph{We are going to construct solutions for GKZ systems with
$\beta \in \MM\,.$ Of special interest for applications to mirror symmetry are
the cases $\beta=0$ and
$\beta=-\ga_0$ with $\ga_0$ as in the definition of reflexive
Gorenstein cone (definition \ref{gorco}).}

The idea is as follows.
Gel'fand-Kapranov-Zelevinskii \cite{gkz1} give
solutions for (\ref{eq:gkzlin})-(\ref{eq:gkzmon}) in the form of so-called
$\Gamma$-series
\begin{equation}\label{eq:gaser}
\sum_{\ell\in\LL}\;\prod_{j=1}^N
\frac{\:v_j^{\gamma_j+\ell_j}}{\Gamma (\gamma_j+\ell_j+1)}
\end{equation}
$\Gamma$ is the usual $\Gamma$-function,
$\ell=(\ell_1,\ldots,\ell_N)^t\in\LL\subset\ZZ^N$.
The series depends on  an additional parameter
$\gamma=(\gamma_1,\ldots,\gamma_N)^t\in\CC^N$ which must satisfy
\begin{equation} \label{eq:Agabe}
\gamma_1\ga_1+\ldots+\gamma_N\ga_N\,=\,\beta
\end{equation}
Allowing the obvious formal rules for differentiating such $\Gamma$-series one
sees that
the functional equations of the $\Gamma$-function guarantee that
(\ref{eq:gaser}) satisfies the differential equations (\ref{eq:gkzmon})
and that
condition (\ref{eq:Agabe}) on $\gamma$ takes care of (\ref{eq:gkzlin}).
\emph{The issue is to interpret the $\Gamma$-series (\ref{eq:gaser}) as a
function on some domain.}
In order that (\ref{eq:gaser}) can be realized as a function $\gamma$ must
satisfy more conditions.
Gel'fand-Kapranov-Zelevinskii obtain convenient conditions from
a triangulation $\gT$ of the convex hull of
$\{\ga_1,\ldots,\ga_N\}$. However, if $\beta$ is in $\MM$ and the triangulation
has more than one maximal simplex,
the vectors $\gamma$ which satisfy these extra conditions do not
provide enough $\Gamma$-series solutions for the GKZ system. This
phenomenon is called \emph{resonance} \cite{gkz1}. An extreme
case of resonance, in which all $\Gamma$-series coincide,
occurs when $\beta$ is in $\MM$ and $\gT$ is unimodular.

\begin{definition} (cf. {} \cite{stur})
A triangulation is called \emph{unimodular}
if all its maximal simplices have volume $1\,;$ the volume of a maximal simplex
$\conv\{\ga_{i_1},\ldots,\ga_{i_\na}\}$ is defined as
$|\det(\ga_{i_1},\ldots,\ga_{i_\na})|\,.$
\end{definition}

\

To get around the resonance problem for $\beta \in \MM$ we proceed as follows.
Fixing a solution $\gamma^\circ\in\ZZ^N$ for equation (\ref{eq:Agabe}) we write
the general solution of (\ref{eq:Agabe}) as
$\gamma=\gamma^\circ+\ec$ with $\ec=(c_1,\ldots,c_N)^t$ such that
\begin{equation} \label{eq:linc0}
c_1\ga_1+\ldots+c_N\ga_N=0
\end{equation}
and note $\gamma+\LL=\ec+\cA^{-1}(\beta)\,.$ Thus (\ref{eq:gaser})
becomes
$\sum_{\lambda\in\cA^{-1}(\beta)}\hspace{2mm}\prod_{j=1}^N
\frac{\:v_j^{c_j+\lambda_j}}{\Gamma (c_j+\lambda_j+1)}\:.$
Multiplying this by $\prod_{j=1}^N\Gamma (c_j+1)$ we obtain
\begin{equation}\label{eq:form1}
\Phi_{\gT,\beta}(\ev)\::=\:\sum_{\lambda\in\cA^{-1}(\beta)}\hspace{2mm}
Q_\lambda(\ec)\;\cdot\;\prod_{j=1}^N v_j^{\lambda_j}\;\cdot\;\prod_{j=1}^N
v_j^{c_j}
\end{equation}
where
\begin{equation}\label{eq:bico1}
Q_\lambda(\ec)\::=\:
\frac{\prod_{\lambda_j<0}\prod_{k=0}^{-\lambda_j-1}(c_j-k)}
{\prod_{\lambda_j>0}\prod_{k=1}^{\lambda_j}(c_j+k)}\,.
\end{equation}
The key observation is that (\ref{eq:bico1}) and (\ref{eq:form1}) also make
sense when $c_1,\ldots,c_N$ are taken from a $\QQ$-algebra in which they are
nilpotent. The expression $v_j^{c_j}$ can still be interpreted as $\exp(c_j\log
v_j)$.

\

\begin{definition}\label{ring}
Let $\sA=(a_{ij})$ denote the $\na\times N$-matrix with columns
$\ga_1,\ldots,\ga_N$. For a regular triangulation $\gT$
(cf. $\S$ \ref{triangs})
of the polytope $\Delta\,:=\conv\{\ga_1,\ldots,\ga_N\}$ we define:
\begin{equation}\label{eq:ring}
\cR_{\sA,\gT}\::=\:\modquot{\ZZ[D^{-1}][C_1,\ldots,C_N]}{\cJ}
\end{equation}
where $\cJ$ is the ideal generated by the linear forms
\begin{equation} \label{eq:linrel}
a_{i1}C_1+\ldots+a_{iN}C_N\hspace{3mm}\textrm{ for } i=1,\ldots,\na
\end{equation}
and by the monomials
\begin{equation}\label{eq:srrel}
C_{i_1}\cdot\ldots\cdot C_{i_s} \hspace{4mm}\textrm{ with }
 \hspace{4mm} \conv\{\ga_{i_1},\ldots,\ga_{i_s}\}
 \hspace{4mm}\textrm{ not a simplex in }\gT\, ;
\end{equation}
$D$ is the product of the volumes of the maximal simplices of $\gT\,.$

We write $c_i$ for the image of $C_i$ in $\cR_{\sA,\gT}\,.$
\end{definition}

\

In theorem \ref{dims} we show that
$\cR_{\sA,\gT}$ is a free $\ZZ[D^{-1}]$-module with rank equal to the number of
maximal simplices in the triangulation. This implies that $c_1,\ldots,c_N$ are
nilpotent and hence
$$
Q_\lambda(\ec)\in\cR_{\sA,\gT}\otimes\QQ
$$
\begin{theorem}\label{values}\nonumber
With this interpretation of $Q_\lambda(\ec)$
the function $\Phi_{\gT,\beta}(\ev)$\\ defined by (\ref{eq:form1}) takes values
in the
\emph{ algebra } $\cR_{\sA,\gT}\otimes\CC\,.$\qed
\end{theorem}

\

The domain of definition of the function $\Phi_{\gT,\beta}(\ev)$
is discussed hereafter.
Relation (\ref{eq:linrel}) ensures that this function
$\Phi_{\gT,\beta}(\ev)$ satisfies
the differential equations (\ref{eq:gkzlin}); it automatically satisfies
(\ref{eq:gkzmon}). Relation (\ref{eq:srrel}) ensures that
the series expansion for $\Phi_{\gT,\beta}(\ev)$ only contains $\lambda$'s
(i.e. $Q_\lambda(\ec)\neq 0$) which satisfy
\begin{equation}\label{eq:support}
\cA\,\lambda=\beta\hspace{4mm}\textrm{and}\hspace{4mm}
\conv\{\ga_i\:|\:\lambda_i<0\}\hspace{3mm}
\textrm{is a simplex in the triangulation  } \gT\,.
\end{equation}
This is important for determining a domain of definition for
$\Phi_{\gT,\beta}(\ev)$.

As we tried to distinguish a kind of regular behavior for the
$\lambda$'s which satisfy (\ref{eq:support}), we were led to triangulations
for which the intersection of the maximal simplices is not empty. We call
\begin{equation}\label{eq:introcore}
\core\gT\::=\textrm{ intersection of the maximal simplices of }\gT
\end{equation}
the core of the triangulation $\gT$.
We use the short notation  $i\in\core\gT$ for $\ga_i\in\core\gT$.
The following result is corollary \ref{corid} in section \ref{coresection}.

\

\begin{theorem}\label{hascore}
Assume $\core\gT\neq\emptyset$ and $\beta=\sum_{i\in\core\gT} m_i\ga_i$
with all $m_i<0$. Then the function $\Phi_{\gT,\beta}(\ev)$
takes values in the
\emph{ principal ideal } $c_{\core}\cR_{\sA,\gT}\otimes\CC$
where
$$
c_{\core}\::=\:\prod_{i\in\core\gT} c_i
$$
Multiplication by $c_{\core}$ on $\cR_{\sA,\gT}$ induces a linear
isomorphism
\begin{equation}\label{eq:coriso}
\modquot{\cR_{\sA,\gT}}{\Ann c_{\core}}
\stackrel{\simeq}{\longrightarrow}c_{\core}\cR_{\sA,\gT}
\end{equation}
Thus one can also say that the function $\Phi_{\gT,\beta}(\ev)$
takes values in the algebra
$\modquot{\cR_{\sA,\gT}}{\Ann c_{\core}}\otimes\CC\,.$

\qed
\end{theorem}

\

By composing $\Phi_{\gT,\beta}$ with a linear map $\cR_{\sA,\gT}\rightarrow\CC$
one obtains a $\CC$-multi-valued function which satisfies the system of
differential equations (\ref{eq:gkzlin})-(\ref{eq:gkzmon}).
When $\beta=0$ and  $\gT$ is unimodular all solutions of
(\ref{eq:gkzlin})-(\ref{eq:gkzmon}) can be obtained in this way; see theorem
\ref{isosol}.

For $\beta\neq 0$ not all solutions of
(\ref{eq:gkzlin})-(\ref{eq:gkzmon}) can be
obtained in this way. Yet what we need for mirror symmetry are the
solutions which can be obtained in this way for appropriate $\beta$ and $\gT$;
see theorem \ref{mainthm}. Our proof of this theorem makes essential use
of the relation:
\begin{equation}\label{eq:recur}
\frac{\partial}{\partial v_i}\Phi_{\gT,\beta}(\ev)\:=\:
\Phi_{\gT,\beta-\ga_i}(\ev)\,.
\end{equation}
which follows imediately from the formulas (\ref{eq:form1}) and
(\ref{eq:bico1}).

\begin{remark}\textup{
The ideal generated by the monomials in (\ref{eq:srrel}) is known as the
\emph{Stanley-Reisner ideal} and has been defined
for finite simplicial complexes in general \cite{stan}.
It is well-known \cite{dani,ful,oda} that the cohomology ring of a toric
variety constructed from a complete simplicial fan has a presentation by
generators and relations as in (\ref{eq:linrel})-(\ref{eq:srrel}).
Unimodular triangulations whose core is not empty and is
not contained in the boundary of $\Delta$,
give rise to such toric varieties and in
that case $\cR_{\sA,\gT}$ is indeed the
cohomology ring of a toric variety;
see theorem \ref{protor}. However not all triangulations to which the present
discussion applies are of this kind. For instance for the triangulation
$\gT_5$ in figure 1 we find $\cR_{\sA,\gT_5}=
\ZZ[c_1,c_2,c_5]/(c_1^2,\,c_2^2,\,c_5^2,\,c_1c_2,\,c_2c_5)$.
An element like $c_2$ which annihilates the whole degree $1$ part of
$\cR_{\sA,\gT_5}$ can not exist in the cohomology of a toric variety.}
\end{remark}

\begin{remark}\textup{
Our method for solving GKZ systems in the resonant case evolved directly from
the $\Gamma$-series of Gel'fand-Kapranov-Zelevinskii.
In hindsight it can also be viewed as a variation on the classical
method of Frobenius  \cite{frob}. The latter would
view $\gamma_1,\ldots,\gamma_N$ in (\ref{eq:gaser}) or $c_1,\ldots,c_N$
in (\ref{eq:form1}) as variables with a restriction given by (\ref{eq:Agabe})
or (\ref{eq:linc0}); then differentiate (repeatedly if necessary) with respect
to these variables and set $\gamma=(\gamma_1,\ldots,\gamma_N)$ in the
derivatives equal to its special value $\gamma^\circ$, c.q. set
$c_1=\ldots=c_N=0$, to obtain solutions for
(\ref{eq:gkzlin})-(\ref{eq:gkzmon}).
Frobenius \cite{frob} considered only functions in one variable.
In the case with more variables one also needs a good bookkeeping device for
the linear relations between the solutions of the differential equations. The
rings $\cR_{\sA,\gT}$ resp. $\modquot{\cR_{\sA,\gT}}{\Ann c_{\core}}$
are such a bookkeeping devices.
Hosono-Klemm-Theisen-Yau have applied Frobenius' method
directly in the situation of the Picard-Fuchs equations of certain families of
Calabi-Yau threefolds; see   \cite{HKTY} formulas (4.9)
and (4.10). In their work the cohomology ring of the mirror Calabi-Yau
threefold plays a similar role of bookkeeper; in fact
$\modquot{\cR_{\sA,\gT}}{\Ann c_{\core}}$ is the cohomology ring of the mirror
Calabi-Yau manifold . The way in which we arrive at our result looks quite
different from that in  \cite{HKTY}
$\S 4$. Moreover the formulation in op. cit. is restricted to the situation of
Calabi-Yau threefolds.}
\end{remark}

\begin{remark}\textup{
Some of our $\Phi_{\gT,\beta}$'s are similar to expressions
presented by Givental in  \cite{giv} theorems 3 and 4; more
specifically,  $\vec{g}_l$ in \cite{giv} thm. 4
is a special case of $\Phi_{\gT,\beta}$ in our theorem
\ref{hascore} with $\beta=-\sum_{i\in\core\gT}\ga_i$,
whereas in \cite{giv} thm. 3 there is a difference in that
the input data
are not subject to (\ref{eq:hyperplane}) in condition \ref{gkzcond}.
The algebra $H$ in \cite{giv} thm. 3 is the cohomology algebra of a
toric variety while our $\cR_{\sA,\gT}$ for
appropriate $\gT$ is also the cohomology algebra of a toric variety.
The algebra $H$ in \cite{giv} thm. 4 is
the algebra $\modquot{\cR_{\sA,\gT}}{\Ann c_{\core}}$
in our theorem \ref{hascore}.
}
\end{remark}
$$
$$
For a proper treatment of the logarithms which appear in (\ref{eq:form1}) we
set
\begin{equation}\label{eq:logcoord}
\begin{array}{lcl}
v_j&:=&\exp(\tp z_j) \hspace{4mm} (j=1,\ldots,N)\\
\ez&:=&(z_1,\ldots,z_N)\:\in\: \CC^{N\vee}\\
\ec&:=&\,(c_1,\ldots,c_N)^t\in\cR_{\sA,\gT}\,\otimes\ZZ^N\,;
\end{array}
\end{equation}
by (\ref{eq:linc0}) $\; \ec$ lies
in fact in $\cR_{\sA,\gT}\,\otimes\LL\,.$
Instead of (\ref{eq:form1}) we now consider
\begin{equation}\label{eq:form2}
\Psi_{\gT,\beta}(\ez)\::=\:
\sum_{\lambda\in\cA^{-1}(\beta)}\hspace{2mm}
Q_\lambda(\ec)\ee^{\tp\ez\cdot\lambda}\cdot\ee^{\tp\ez\cdot\ec}\,.
\end{equation}
Note that $\ee^{\tp\ez\cdot\ec}$ is just a polynomial, but
$\sum_{\lambda\in\cA^{-1}(\beta)}\;
Q_\lambda(\ec)\ee^{\tp\ez\cdot\lambda}$ is really a series. In section
\ref{solutions} we analyse the convergence of this series and give a
domain $\cV_\gT$ in $\CC^{N\vee}$ on which the function
$\Psi_{\gT,\beta}$ is defined; see theorem \ref{defdom}.

The domain $\cV_\gT$ is invariant under translations by elements of
$\ZZ^{N\vee}$ and by elements
of $\MM_\CC^\vee\::=\:\Hom{\MM}{\CC}\subset\CC^{N\vee}$.
{}From (\ref{eq:form2}) one immediately sees
\begin{eqnarray}\label{eq:monodromy}
\Psi_{\gT,\beta}(\ez+\mu)&=&\ee^{\tp\mu\,\cdot\ec}\cdot
\Psi_{\gT,\beta}(\ez)\hspace{5mm}\forall \:\mu\in\ZZ^{N\vee}
\\
\label{eq:character}
\Psi_{\gT,\beta}(\ez+\sfm)&=&\ee^{\tp\sfm\,\cdot\beta}\cdot
\Psi_{\gT,\beta}(\ez)
\hspace{5mm}\forall \:\sfm\in\MM_\CC^\vee
\,.
\end{eqnarray}
The functional equation (\ref{eq:monodromy}) gives the monodromy of
$\Phi_{\gT,\beta}$, when viewed as a multivalued function on
$\modquot{\cV_\gT}{\ZZ^{N\vee}}$ with values
in the vector space $\cR_{\sA,\gT}\otimes\CC$.
Because of (\ref{eq:linrel}) elements of $\MM_\ZZ^\vee\::=\:\Hom{\MM}{\ZZ}$
give trivial monodromy and the actual monodromy comes from
$\LL_\ZZ^\vee\::=\:\Hom{\LL}{\ZZ}$.

As $\MM_\ZZ^\vee$ acts trivially, the translation action of
$\MM_\CC^\vee$ descends to an action of the torus
$\modquot{\MM_\CC^\vee}{\MM_\ZZ^\vee}=\Hom{\MM}{\CC^\ast}$. The
functional equation (\ref{eq:character}), whose infinitesimal analogues
are the differential equations (\ref{eq:gkzlin}),
means that $\Psi_{\gT,\beta}$ is an eigenfunction with
character $\beta$.

If one wants an invariant function for $\beta\neq 0$
one must replace the range of values of $\Psi_{\gT,\beta}$ by
$\modquot{(\cR_{\sA,\gT}\otimes\CC)}{\CC^{\,\ast}}$, the orbit space
for the natural $\CC^{\,\ast}$-action on the vector space
$\cR_{\sA,\gT}\otimes\CC$. On a possibly slightly smaller domain of
definition the invariant function even takes values in the
projective space $\PP(\cR_{\sA,\gT}\otimes\CC)$.
The $\MM_\CC^\vee$-invariant function $\Psi_{\gT,\beta}\bmod\CC^{\,\ast}$
is defined on the domain
$\LL_\RR^\vee+\sqrt{-1} \cB_\gT$ in $\LL_\CC^\vee$; cf. formula
(\ref{eq:domain4}). The (multivalued)
function $\Phi_{\gT,\beta}\bmod\CC^{\,\ast}$ is defined on a domain in the
torus $\Hom{\LL}{\CC^{\,\ast}}$.

For a good overall picture it is appropriate to point out here that the
pointed secondary fan (the construction of which is recalled in section
\ref{secfan}) defines a toric variety which compactifies the torus
$\Hom{\LL}{\CC^{\,\ast}}$. To each regular triangulation of $\Delta$
corresponds a special point in the boundary of this compactification.
The domain of definition of $\Phi_{\gT,\beta}\bmod\CC^{\,\ast}$
is the intersection of the torus $\Hom{\LL}{\CC^{\,\ast}}$ and a neighborhood
of the special point corresponding to $\gT$;
see the end of section \ref{solutions}.

\begin{example}\textup{
Let $\ga_1,\ldots,\ga_6$ be the columns of the following matrix $\sA$:}
\end{example}

\begin{picture}(300,80)(-90,0)
\put(-75,35){\makebox(100,50){$\sA=\left(
\begin{array}{rrrrrr}
1&1&1&1&1&1\\
0&1&-1&0&1&0\\
1&1&0&0&0&-1
\end{array}
\right)$}}
\put(-75,-5){\makebox(100,50){$\sB=\left(
\begin{array}{rrrrrr}
1&0&0&-2&0&1\\
0&1&1&-3&0&1\\
0&0&1&-2&1&0
\end{array}
\right)$}}
\put(70,0){\makebox(50,70){$\Delta=$}}
\put(145,10){\line(1,1){25}}
\put(145,10){\line(-1,1){25}}
\put(145,60){\line(1,0){25}}
\put(145,60){\line(-1,-1){25}}
\put(170,35){\line(0,1){25}}
\put(145,10){\circle*{5}}
\put(120,35){\circle*{5}}
\put(145,35){\circle*{5}}
\put(170,35){\circle*{5}}
\put(145,60){\circle*{5}}
\put(170,60){\circle*{5}}
\put(133,53){\makebox(10,20){$\ga_1$}}
\put(172,53){\makebox(10,20){$\ga_2$}}
\put(115,33){\makebox(10,20){$\ga_3$}}
\put(143,33){\makebox(10,20){$\ga_4$}}
\put(172,31){\makebox(10,20){$\ga_5$}}
\put(149,0){\makebox(10,20){$\ga_6$}}
\end{picture}

\noindent These satisfy conditions (\ref{eq:aspan}) and
(\ref{eq:hyperplane}) with $\MM=\ZZ^3$ and
$\ga_0^\vee=(1,0,0)$.

Figure 1 shows all regular triangulations of the polytope $\Delta$,
with two triangulations joined by an edge iff the corresponding
cones in the pointed secondary fan are adjacent.
\setlength{\unitlength}{.7pt}

\begin{picture}(250,340)(-50,-40)
\put(50,-30){\makebox(160,35){\textbf{figure 1}}}
\put(-15,-20){\makebox(60,20){$\gT_4$}}
\put(185,-20){\makebox(60,20){$\gT_1$}}
\put(265,20){\makebox(60,20){$\gT_2$}}
\put(-15,229){\makebox(60,20){$\gT_9$}}
\put(75,269){\makebox(60,20){$\gT_8$}}
\put(265,269){\makebox(60,20){$\gT_7$}}
\put(90,65){\makebox(60,20){$\gT_3$}}
\put(95,182){\makebox(60,20){$\gT_{10}$}}
\put(210,125){\makebox(60,20){$\gT_5$}}
\put(200,228){\makebox(60,20){$\gT_6$}}

\thicklines
\put(15,34){\line(0,1){160}}
\put(33,15){\line(1,0){162}}
\put(36,23){\line(2,1){15}}
\put(64,37){\line(2,1){15}}

\put(34,215){\line(1,0){82}}
\put(32,223){\line(2,1){51}}
\put(215,34){\line(0,1){83}}
\put(232,23){\line(2,1){51}}
\put(270,240){\line(2,1){15}}
\multiput(95,75)(0,24){7}{\line(0,1){12}}
\multiput(117,55)(25,0){4}{\line(1,0){12}}
\multiput(218,55)(25,0){2}{\line(1,0){12}}
\put(266,55){\line(1,0){10}}
\put(277,255){\line(-1,0){162}}
\put(295,237){\line(0,-1){162}}
\put(154,215){\line(6,1){88}}
\put(218,153){\line(2,5){29}}
\put(205,145){\line(-1,1){60}}

\thinlines
\put(15,0){\line(1,1){15}}
\put(15,0){\line(-1,1){15}}
\put(15,30){\line(1,0){15}}
\put(15,30){\line(-1,-1){15}}
\put(30,15){\line(0,1){15}}
\put(15,0){\circle*{3}}
\put(0,15){\circle*{3}}
\put(15,15){\circle*{3}}
\put(30,15){\circle*{3}}
\put(15,30){\circle*{3}}
\put(30,30){\circle*{3}}
\put(15,15){\line(-1,0){15}}
\put(15,15){\line(1,1){15}}
\put(15,0){\line(0,1){30}}
\put(15,0){\line(1,2){15}}

\put(15,200){\line(1,1){15}}
\put(15,200){\line(-1,1){15}}
\put(15,230){\line(1,0){15}}
\put(15,230){\line(-1,-1){15}}
\put(30,215){\line(0,1){15}}
\put(15,200){\circle*{3}}
\put(0,215){\circle*{3}}
\put(15,215){\circle*{3}}
\put(30,215){\circle*{3}}
\put(15,230){\circle*{3}}
\put(30,230){\circle*{3}}
\put(15,200){\line(0,1){30}}
\put(15,200){\line(1,2){15}}

\put(215,0){\line(1,1){15}}
\put(215,0){\line(-1,1){15}}
\put(215,30){\line(1,0){15}}
\put(215,30){\line(-1,-1){15}}
\put(230,15){\line(0,1){15}}
\put(215,0){\circle*{3}}
\put(200,15){\circle*{3}}
\put(215,15){\circle*{3}}
\put(230,15){\circle*{3}}
\put(215,30){\circle*{3}}
\put(230,30){\circle*{3}}
\put(200,15){\line(1,0){30}}
\put(215,0){\line(0,1){30}}
\put(215,15){\line(1,1){15}}

\put(95,40){\line(1,1){15}}
\put(95,40){\line(-1,1){15}}
\put(95,70){\line(1,0){15}}
\put(95,70){\line(-1,-1){15}}
\put(110,55){\line(0,1){15}}
\put(95,40){\circle*{3}}
\put(80,55){\circle*{3}}
\put(95,55){\circle*{3}}
\put(110,55){\circle*{3}}
\put(95,70){\circle*{3}}
\put(110,70){\circle*{3}}
\put(95,55){\line(-1,0){15}}
\put(95,55){\line(0,-1){15}}
\put(95,55){\line(1,1){15}}
\put(110,70){\line(-2,-1){30}}
\put(110,70){\line(-1,-2){15}}

\put(295,40){\line(1,1){15}}
\put(295,40){\line(-1,1){15}}
\put(295,70){\line(1,0){15}}
\put(295,70){\line(-1,-1){15}}
\put(310,55){\line(0,1){15}}
\put(295,40){\circle*{3}}
\put(280,55){\circle*{3}}
\put(295,55){\circle*{3}}
\put(310,55){\circle*{3}}
\put(295,70){\circle*{3}}
\put(310,70){\circle*{3}}
\put(295,55){\line(0,-1){15}}
\put(295,55){\line(1,1){15}}
\put(280,55){\line(1,0){30}}
\put(280,55){\line(2,1){30}}

\put(95,240){\line(1,1){15}}
\put(95,240){\line(-1,1){15}}
\put(95,270){\line(1,0){15}}
\put(95,270){\line(-1,-1){15}}
\put(110,255){\line(0,1){15}}
\put(95,240){\circle*{3}}
\put(80,255){\circle*{3}}
\put(95,255){\circle*{3}}
\put(110,255){\circle*{3}}
\put(95,270){\circle*{3}}
\put(110,270){\circle*{3}}
\put(110,270){\line(-1,-2){15}}
\put(110,270){\line(-2,-1){30}}

\put(295,240){\line(1,1){15}}
\put(295,240){\line(-1,1){15}}
\put(295,270){\line(1,0){15}}
\put(295,270){\line(-1,-1){15}}
\put(310,255){\line(0,1){15}}
\put(295,240){\circle*{3}}
\put(280,255){\circle*{3}}
\put(295,255){\circle*{3}}
\put(310,255){\circle*{3}}
\put(295,270){\circle*{3}}
\put(310,270){\circle*{3}}
\put(280,255){\line(1,0){30}}
\put(280,255){\line(2,1){30}}

\put(135,200){\line(1,1){15}}
\put(135,200){\line(-1,1){15}}
\put(135,230){\line(1,0){15}}
\put(135,230){\line(-1,-1){15}}
\put(150,215){\line(0,1){15}}
\put(135,200){\circle*{3}}
\put(120,215){\circle*{3}}
\put(135,215){\circle*{3}}
\put(150,215){\circle*{3}}
\put(135,230){\circle*{3}}
\put(150,230){\circle*{3}}
\put(135,230){\line(0,-1){30}}
\put(135,230){\line(1,-1){15}}

\put(215,120){\line(1,1){15}}
\put(215,120){\line(-1,1){15}}
\put(215,150){\line(1,0){15}}
\put(215,150){\line(-1,-1){15}}
\put(230,135){\line(0,1){15}}
\put(215,120){\circle*{3}}
\put(200,135){\circle*{3}}
\put(215,135){\circle*{3}}
\put(230,135){\circle*{3}}
\put(215,150){\circle*{3}}
\put(230,150){\circle*{3}}
\put(200,135){\line(1,0){30}}
\put(215,120){\line(0,1){30}}
\put(230,135){\line(-1,1){15}}

\put(255,220){\line(1,1){15}}
\put(255,220){\line(-1,1){15}}
\put(255,250){\line(1,0){15}}
\put(255,250){\line(-1,-1){15}}
\put(270,235){\line(0,1){15}}
\put(255,220){\circle*{3}}
\put(240,235){\circle*{3}}
\put(255,235){\circle*{3}}
\put(270,235){\circle*{3}}
\put(255,250){\circle*{3}}
\put(270,250){\circle*{3}}
\put(270,235){\line(-1,0){30}}
\put(270,235){\line(-1,1){15}}
\end{picture}

The columns of the matrix $\sB^t$
constitute a $\ZZ$-basis for $\LL$ by means of which one can identify
$\LL$ with $\ZZ^3$ and $\LL_\RR^\vee$ with  $\RR^{3\vee}$.
The rows $\gb_1,\ldots,\gb_6$ of the matrix $\sB^t$ are then
identified with the images of the standard basis vectors
under the projection $\RR^{6\vee}\rightarrow\LL_\RR^\vee$, dual
to the inclusion $\LL\subset\ZZ^6$.
Thus one finds $\ell_j=\gb_j\cdot\ell$ for every
$\ell\in\LL_\RR\simeq\RR^3$ and (\ref{eq:support}) becomes a condition
on the signs of $\gb_1\cdot\ell+\gamma^\circ_1,\ldots,
\gb_6\cdot\ell+\gamma^\circ_6\,.$ The signs give a vector in
$\,\{-1,0,+1\}^6\,$

These sign vectors correspond exactly to the various
strata in the stratification of $\RR^3$ induced by the six
planes $\gb_j\cdot\ex+\gamma^\circ_j=0$  ($j=1,\ldots,6$).
Figure 2 shows the zonotope spanned by $\gb_1,\ldots,\gb_6$.
The $3-j$-dimensional faces of this zonotope correspond bijectively
with the $j$-dimensional strata in the stratification for
$\gamma^\circ=0$. The stratum with sign vector
$(s_1,\ldots,s_6)$ corresponds with the face whose centre is
$s_1\gb_1+s_2\gb_2+s_3\gb_3+s_4\gb_4+s_5\gb_5+s_6\gb_6\,.$
The vertices $1$--$14$ (resp. $15$--$28$)
of the zonotope have sign vectors $(s_1,\ldots,s_6)$
(resp. $-(s_1,\ldots,s_6)$) as given in table 1.

The sign vectors of all faces of the zonotope give all possible
signs for $\ell=(\ell_1,\ldots,\ell_6)\in\LL\,.$ Thus by comparing
this with (\ref{eq:support}) one can see for every triangulation $\gT$
what types of terms are involved in the series of $\Psi_{\gT,0}$.
For example for triangulation $\gT_1$ the series of $\Psi_{\gT_1\!,0}$
involves precisely those $\ell\in\LL$ whose sign vector
corresponds to a face of the zonotope containing at least
one of the vertices $1$, $2$, $3$ or $4$.

\begin{picture}(250,250)(-50,-270)
\put(50,-240){\makebox(250,250){$
\begin{array}{cc}
\epsfig{file=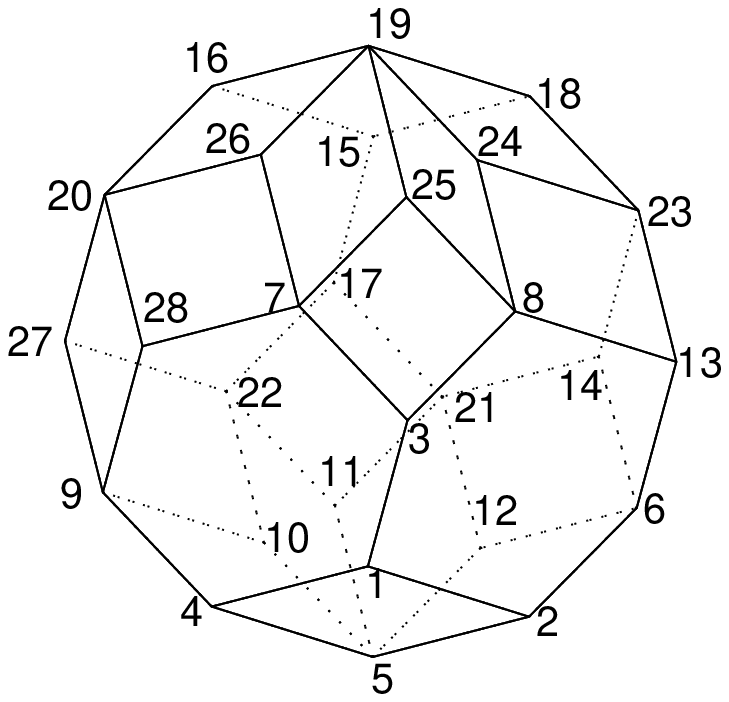,height=61mm} &
\begin{array}{rrrrrrrl}
     1 &+ &+ &+ &- &+ &+ & 15\\
     2 &- &+ &+ &- &+ &+ & 16\\
     3 &+ &- &+ &- &+ &+ & 17\\
     4 &+ &+ &+ &- &- &+ & 18\\
     5 &- &+ &+ &- &- &+ & 19\\
     6 &- &+ &+ &- &+ &- & 20\\
     7 &+ &- &- &- &+ &+ & 21\\
     8 &+ &- &+ &- &+ &- & 22\\
     9 &+ &+ &- &- &- &+ & 23\\
    10 &- &+ &- &- &- &+ & 24\\
    11 &- &+ &+ &+ &- &+ & 25\\
    12 &- &+ &+ &- &- &- & 26\\
    13 &- &- &+ &- &+ &- & 27\\
    14 &- &+ &+ &+ &+ &- & 28
\end{array}
\\ & \\
\textbf{figure 2}&\textbf{table 1}
\end{array}
$}}
\end{picture}
$$ $$
The series
$\Psi_{\gT_1\!,-\ga_4}$ involves the same $\ell$'s with exception of
$\ell=0$ (which corresponds to the $3$-dimensional the zonotope itself).
Using the Pochhammer symbol notation
$(x)_m:=x\,(x+1)\cdot\ldots\cdot(x+m-1)$ we have
\begin{eqnarray*}
\Psi_{\gT_1\!,-\ga_4}&=& c_4\:\ee^{-2\pi i z_4}\:
T_1^{c_1}T_2^{c_2}T_5^{c_5}\times\\
&\times&\hspace{-2mm}\left\{\hspace{2mm}
\sum_{p,q,r\geq 0}(-1)^q
\frac{(2c_1+3c_2+2c_5+1)_{2p+3q+2r}}{(c_1)_p(c_2)_q
(c_5)_r(c_1+c_2)_{p+q}(c_2+c_5)_{q+r}} T_1^pT_2^qT_5^r
\right.\\
&&\hspace{5mm}-c_1\sum_{r\geq 0,-q\leq p<0}(-1)^{q+p}
\frac{(2p+3q+2r)!(-p-1)!}{q!r!(p+q)!(q+r)!}
 T_1^{p}T_2^{q}T_5^{r}
\\
&&\hspace{5mm}-c_5\sum_{p\geq 0,-q\leq r<0}(-1)^{q+r}
\frac{(2p+3q+2r)!(-r-1)!}{p!q!(p+q)!(q+r)!}
 T_1^{p}T_2^{q}T_5^{r}
\\
&&\left.\hspace{5mm}-c_2
\sum_{{-p\leq q<0\,,-r\leq q<0}}
\frac{(2p+3q+2r)!(-q-1)!}{p!r!(p+q)!(q+r)!}
 T_1^{p}T_2^{q}T_5^{r}\hspace{2mm} \right\}
\end{eqnarray*}
where
$$
T_1:=\ee^{2\pi i(z_1-2z_4+z_6)}\;,\hspace{3mm}
T_2:=\ee^{2\pi i(z_2+z_3-3z_4+z_6)}\;,\hspace{3mm}
T_5:=\ee^{2\pi i(z_3-2z_4+z_5)}
$$
$$
c_4=-2c_1-3c_2-2c_5
$$
and
$$
\cR_{\sA,\gT_1}=\modquot{\ZZ[c_1,c_2,c_5]}{(c_1^2-c_2^2,\,
c_1^2-c_5^2,\,c_1^2+c_1c_2,\,c_1^2+c_2c_5,\,c_1c_5)}\,.
$$
Note that $c_4c_1=c_4c_2=c_4c_5$. One may therefore simplify
the expression for $\Psi_{\gT_1\!,-\ga_4}$ and replace
$c_2$ and $c_5$ by $c_1$.


\section{Regular triangulations and the pointed secondary fan}
\label{regtrifan}

In this section we review some results about regular triangulations
and about the pointed secondary fan, essentially following \cite{bfs}.
One may take as a definition of {\it regular triangulations} that these are the
triangulations produced by the construction in this section; see in particular
proposition \ref{protri}.

\subsection{Regular triangulations}
\label{triangs}

We start from  a set of vectors
$\{\ga_1,\ldots,\ga_N\}$ in $\ZZ^\na$ satisfying condition
\ref{gkzcond}. Let $\Delta=\conv\{\ga_1,\ldots,\ga_N\}$ denote
the convex hull of this set of points in $\RR^\na$.
We are interested in triangulations of $\Delta$ such that
all vertices are among  the marked points $\ga_1,\ldots,\ga_N$. The
notation can be conveniently simplified by referring to a simplex
$\conv\{\ga_{i_1},\ldots,\ga_{i_m}\}$
by just the index set $\{i_1,\ldots,i_m\}$. We will allways take
the indices in
increasing order. If $\gT$ is a triangulation,
we write $\gT^m$ for the set of simplices with $m$ vertices.
A triangulation is
completely determined by its set of maximal simplices $\gT^\na$.

For the construction of a regular triangulation we take an $N$-tuple of
positive real numbers $\ed=(d_1,\ldots,d_N)$ and consider the polytope
\begin{equation}\label{eq:pd}
\cP_\ed:=\conv\{0,d_1^{-1}\ga_1,\ldots,d_N^{-1}\ga_N\}\:\subset\:\RR^\na\,.
\end{equation}
Consider a subset $I=\{i_1,\ldots,i_{\na}\}$ of $\{1,\ldots,N\}$ for which
$\ga_{i_1},\ldots,\ga_{i_{\na}}$ are linearly independent.
The affine hyperplane through
$d_{i_1}^{-1}\ga_{i_1},\ldots,d_{i_{\na}}^{-1}\ga_{i_{\na}}$ is given by the
equation $D_{\ed,I}(\ex)=0$ with
\begin{equation}\label{eq:faceq}
D_{\ed,I}(\ex)\::=\:
\det\left(\begin{array}{cccc}
d_{i_1}^{-1}\ga_{i_1}&\ldots&d_{i_{\na}}^{-1}\ga_{i_{\na}}&\ex\\
1&\ldots&1&1\end{array}\right)
\end{equation}
Write $I^\ast:=\{1,\ldots,N\}\setminus I$.
Then $\{d_{i_1}^{-1}\ga_{i_1},\ldots,d_{i_{\na}}^{-1}\ga_{i_{\na}}\}$
lies in a codimension $1$ face of $\cP_\ed$ if and only if for all
$j\in I^\ast$:
\begin{equation}\label{eq:face}
D_{\ed,I}(d_j^{-1}\ga_j)\cdot D_{\ed,I}(0)\geq 0
\end{equation}
This face is a simplex with vertices
$d_{i_1}^{-1}\ga_{i_1},\ldots,d_{i_{\na}}^{-1}\ga_{i_{\na}}$  iff
$D_{\ed,I}(d_j^{-1}\ga_j)\neq 0$ for every $j\in I^\ast$.
Thus if $\ed$ does not lie on any hyperplane in $\RR^N$ given by the vanishing
of $D_{\ed,I}(d_j^{-1}\ga_j)$ for some $I$ and $j$ with $j\not\in I$,
then all faces of $\cP_\ed$
opposite to the vertex $0$ are simplicial.

In this case the projection with center $0$ projects the boundary of $\cP_\ed$
onto \emph{a triangulation $\gT$ of $\Delta$.}
The maximal simplices of $\gT$ are those $I=\{i_1,\ldots,i_{\na}\}$
for which $D_{\ed,I}(d_j^{-1}\ga_j)\cdot D_{\ed,I}(0)>0$ holds for every
$j\in I^\ast\,.$

\

Let $\sA=(a_{ij})$ denote the $\na\times N$-matrix with columns
$\ga_1,\ldots,\ga_N$.
The triangulation obviously depends only on $\ed$ modulo the row space of
$\sA\,.$ Let us reformulate the above construction accordingly.

Take $\LL=\ker\sA\subset\ZZ^N$ as in (\ref{eq:L}). Assumption
(\ref{eq:hyperplane}) implies $\ell_1+\ldots+\ell_N\:=\:0$ for every
$\ell=(\ell_1,\ldots,\ell_N)^t\in\LL\,.$
Take an $\bnl\times N$-matrix $\sB$ with entries in $\ZZ$ such that columns of
$\sB^t$ constitute a basis for $\LL\,.$

Let $w\in\RR^\nl$. Then there exists a row vector of positive real numbers
$\ed=(d_1,\ldots,d_N)$ such that $w\:=\:\sB\ed^t$.
Take the matrices
$$
\widetilde{\sA}\::=\;\left(\begin{array}{c|c} \sA & 0\\ \hrulefill&\hrulefill\\
\ed &1\end{array}\right)
\hspace{5mm}\textrm{ and}\hspace{5mm}
\widetilde{\sB}\::=\:\left(\begin{array}{c|c}\sB&-w\end{array}\right)\,.
$$
Denote by $\widetilde{\sA}_K$ (resp.
$\widetilde{\sB}_K$) the submatrix of $\widetilde{\sA}$  (resp.
$\widetilde{\sB}$) composed of the entries with column index in a subset $K$ of
$\{1,\ldots,N+1\}$.
Since ${\rank\,}\widetilde{\sA}\,=\,\na+1\:,$
${\rank\,}\widetilde{\sB}\,=\,\nl$ and
$\widetilde{\sA}\cdot\widetilde{\sB}^t\:=\:0$
there is a non-zero $r\in\QQ$ such that for every
$J\subset\{1,\ldots,N+1\}$ of cardinality $\na+1$
and $J^\prime=\{1,\ldots,N+1\}\setminus J$
$$
\det(\widetilde{\sA}_J)\,=\,
(-1)^{\sum_{j\in J} j}\;r\,\det(\widetilde{\sB}_{J^\prime})
$$
One sees that (\ref{eq:face}) is equivalent to
\begin{equation}\label{eq:om}
(-1)^{\sharp
\{h\in I^\ast\:|\: h> j\}}
\det\left(\left(\begin{array}{c|c}\sB_{I^\ast\setminus\{j\}}
&w\end{array}\right)\right)
\cdot\det\left(\sB_{I^\ast}\right)\;\geq\;0\,;
\end{equation}
here $\sB_{I^\ast}$ resp. $\sB_{I^\ast\setminus\{j\}}$ is the submatrix of
$\sB$ consisting of the entries with column index in $I^\ast$ resp.
$I^\ast\setminus\{j\}\,.$

\emph{Thus the triangulation $\gT$ can also be constructed from (\ref{eq:om}).}


\subsection{The pointed secondary fan}\label{secfan}

For a more intrinsic formulation which does not refer to a choice of a basis
for $\LL$ we consider the $\bnl$-dimensional real vector space
$\LL_\RR^\vee\::=\:\Hom{\LL}{\RR}\,.$
Let $\gb_1,\ldots,\gb_N\in\LL_\RR^\vee$ be the images of the standard basis
vectors of $\RR^{N\vee}$ under the surjection
$\RR^{N\vee}\twoheadrightarrow\LL_\RR^\vee$
dual to the inclusion $\LL\hookrightarrow\ZZ^N$.
Let $\gB$ (resp. $\gD$) be the collection of those subsets $J$ of
$\{1,\ldots,N\}$ of cardinality $\nl$ (resp. $\nl-1$)
for which the vectors $\gb_j$ ($j\in J$) are linearly independent.
For $K=\{k_1,\ldots,k_\nl\}\in\gB$
and $J=\{j_1,\ldots,j_{\nl-1}\}\in\gD$ we write
\begin{eqnarray*}
\cC_K&:=&\{t_1\gb_{k_1}+\ldots+t_\nl\gb_{k_\nl}\in\LL_\RR^\vee\;|\;
 t_1,\ldots,t_\nl\in\RR_{\geq 0}\;\}\\
H_J&:=&\{t_1\gb_{j_1}+\ldots+t_{\nl-1}\gb_{j_{\nl-1}}\in\LL_\RR^\vee\;|\;
 t_1,\ldots,t_{\nl-1}\in\RR\;\}\,.
\end{eqnarray*}
Choosing a basis for $\LL$ as before one can identify $\LL_\RR^\vee$ with
$\RR^{\nl\vee}$ and
$\gb_1,\ldots,\gb_N$ with the rows of matrix $\sB^t$. The inequality
(\ref{eq:om}) becomes equivalent to the statement
$w\in\cC_{I^\ast}\,.$ The condition $D_{\ed,I}(d_j^{-1}\ga_j)\neq 0$ for the
left hand factor in (\ref{eq:face}) becomes equivalent to
$w\not\in H_J$ for $J=\{1,\ldots,N+1\}\setminus(I\cup\{j\})$.

Thus the preceding discussion shows:

\

\begin{proposition}\label{protri}
\textup{ (cf. \cite{bfs} lemma 4.3.)}
For $w\in\LL_\RR^\vee\setminus\bigcup_{J\in\gD}H_J$ the set
$$
\gT^\na\::=\:\{\,I\:|\:I^\ast\in\gB\textrm{  and  }w\in\cC_{I^\ast}\,\}
$$
is the set of maximal simplices of a regular triangulation $\gT$ of
$\Delta\,$.\\
\textup{(Recall the notation
$I^\ast:=\{1,\ldots,N\}\setminus I\,.$)}
\qed
\end{proposition}

\

If $\gT$ is a regular triangulation of $\Delta$ write
\begin{equation}\label{eq:ct}
\cC_{\gT}\:=\:\bigcap_{I\in\gT^\na} \cC_{I^\ast}\,.
\end{equation}
Then every $w\in\cC_{\gT}\setminus\bigcup_{J\in\gD}H_J$ leads by the above
construction to the same triangulation $\gT$.

The cones $\cC_\gT$ one obtains in this way from all regular triangulations of
$\Delta$ constitute the collection of maximal cones of a complete fan in
$\LL_\RR^\vee\,.$
This fan is called \emph{ the pointed secondary fan}.

\

\begin{remark}\label{delzant}\textup{
The dual (or polar) set of $\cP_\ed$ in (\ref{eq:pd}) is
(e.g. \cite{bat1} def.4.1.1,   \cite{ful} p.24)
\begin{equation}\label{eq:dupo}
\cP_\ed^\vee\::=\:\{\;\ey\in\RR^{\na\vee}\:|\;\ey\cdot\ex\geq -1
\;\textrm{ for all }\;\ex\in\cP_\ed\;\}
\end{equation}
It is the intersection of half-spaces given by the inequalities
$$
\ey\cdot\ga_i+d_i\geq 0\hspace{5mm} (i=1,\ldots,N)
$$
$\cP_\ed^\vee$ is an unbounded polyhedron. Its vertices correspond with the
codimension $1$ faces of $\cP_\ed$ which do not contain $0$.}

\textup{
Adding to $\ed$ an element of the row space of matrix $\sA$ amounts to just
a translation of the polyhedron $\cP_\ed^\vee$ in $\RR^{\na\vee}$.
If $\ed$ gives rise to a unimodular triangulation, then $\cP_\ed^\vee$ is an
(unbounded) \emph{Delzant polyhedron}
in the sense of   \cite{guillemin} p.8. Thus, by the constructions in
\cite{guillemin}} a point in the real cone $\cC_{\gT}$ for a unimodular regular
triangulation $\gT$ can be interpreted as a parameter for the symplectic
structure of a toric variety.
In view of formula (\ref{eq:domain4}) this applies in particular to the
imaginary part of the variable $\ez$ in (\ref{eq:form2}).
\end{remark}


\section{The ring $\cR_{\sA,\gT}\,.$}
\label{thering}

\begin{theorem}\label{dims}
Consider the ring $\cR_{\sA,\gT}$ as in definition \ref{ring}.
\begin{enumerate}
\item
$\cR_{\sA,\gT}$ is a free $\ZZ[D^{-1}]$-module of rank $\sharp\,\gT^\na\,.$
\item
$\cR_{\sA,\gT}$ is a graded ring. Let $\cR_{\sA,\gT}^{(k)}$ denote its
homogeneous component of degree $k$.
Then the Poincar\'e series of $\cR_{\sA,\gT}$ is:
$$
\sum_{k\geqslant 0} \left(\rank\,\cR_{\sA,\gT}^{(k)}\right)\,\tau^k\:=\:
\sum_{m=0}^{\na}\sharp\, (\,\gT^m\,)\,\tau^m(1-\tau)^{\na-m}
$$
where $\sharp\, (\,\gT^m\,)\,=$ the number of simplices with $m$ vertices;
$\sharp\, (\,\gT^0\,)=1$ by convention. In particular
$$
\cR_{\sA,\gT}^{(k)}=0\hspace{4mm}\textrm{for}\hspace{3mm} k\geq\na\,.
$$
\item $\{c_I\:|\:I\in\gT^\na\}$ is a $\ZZ[D^{-1}]$-basis for $\cR_{\sA,\gT}\,.$
\hspace{3mm} (cf. formula (\ref{eq:cbas}))

\end{enumerate}
\end{theorem}

The \textbf{proof of theorem \ref{dims}} closely follows the proofs of Danilov
( \cite{dani} $\S$ 10) and Fulton ( \cite{ful} $\S$ 5.2) for the analogous
presentation of the Chow ring of a complete simplicial toric variety. We
include a proof here in order check that it needs no reference to algebraic
cycles and also works when the simplicial complex is homeomorphic to a ball
instead of a sphere as in \cite{dani,ful}.

For the construction of a basis for $\cR_{\sA,\gT}$ we choose a
vector $\xi$ in $\Delta$ which should be linearly independent from every
$\na-1$-tuple of vectors in $\{\ga_1,\ldots,\ga_N\}$.
If $I$ is a maximal simplex, then $\{\ga_i\}_{i\in I}$ is a basis of $\RR^\na$
and $\xi\,=\,\sum_{i\in I} x_i\ga_i$ with all $x_i\neq 0\,.$ We define
\begin{eqnarray} \label{eq:imin}
I^-\:&:=&\:\{\,i\in I\:|\: x_i<0\,\}\:, \\
c_I\:&:=&\:\prod_{i\in I^-}\;c_i \:.\label{eq:cbas}
\end{eqnarray}

\

Let $\gT$ be associated with $\ed=(d_1,\ldots,d_N)$ as in section
\ref{triangs}. For $I\in\gT^\na$ let $p_I$ be the positive real number such
that $p_I\xi$ lies in the affine hyperplane through the points
$d_{i}^{-1}\ga_{i}$ with $i\in I$, i.e.
$D_{\ed,I}(p_I\xi)=0\,.$
We may assume that $\ed$ is chosen such that  $p_{I_1}\neq p_{I_2}$ whenever
$I_1\neq I_2\,;$ indeed,  for $I_1\neq I_2$ the equality $p_{I_1}= p_{I_2}$
amounts to a non-trivial linear equation for $d_1,\ldots,d_N\,.$
As in \cite{ful} we define a total ordering on $\gT^\na$:
\begin{equation}\label{eq:ord}
I_1<I_2 \hspace{4mm}\textrm{ iff }\hspace{4mm} p_{I_1}<p_{I_2}\,.
\end{equation}

\begin{lemma}\label{grot}\textup{ (cf. \cite{ful} p.101($\ast$))}
If $I_1^-\subset I_2$ then
$I_1\leq I_2$ .
\end{lemma}
\textbf{proof:}
By definition of $p_{I_1}$ there exist $s_j\in\RR$ such that
$
p_{I_1}\xi\:=\sum_{j\in I_1} s_j d_j^{-1} \ga_j
$
and
$
1\:=\sum_{j\in I_1} s_j \:.
$
If $I_1\neq I_2$ and $I_1^-\subset I_2$ then $s_j>0$ for every $j\in
I_1\setminus I_2\,.$
Using this and (\ref{eq:face}) for $I_2$ one checks:
$D_{\ed,I_2}(p_{I_1}\xi)\cdot D_{\ed,I_2}(0)>0\,.$
This shows that $0$ and $p_{I_1}\xi$ lie on the same side of the affine
hyperplane through the points $d_i^{-1}\ga_i$ with $i\in I_2\,.$ Hence:
$p_{I_1}<p_{I_2}\,.$
\qed

\begin{lemma}\label{tussen}\textup{ (cf. \cite{ful} p.102)}
Let $J$ be a simplex in $\gT$.
Then: $I^-\:\subset\:J\:\subset\:I$ where
$I:=\min\{I^\prime\in\gT^\na\:|\:J\subset I^\prime\}\,.$
\end{lemma}
\textbf{proof:} The conclusion is clear if $I=J$. So assume $I\neq J$ and take
$i\in I\setminus J$. Then $I\setminus\{i\}$ is a codim $1$ simplex in the
triangulation, which either is contained in the boundary of $\Delta$ or is
contained in another maximal simplex $I^\prime\neq I$.

If $I\setminus\{i\}$ is a boundary simplex, then $\xi$ and $\ga_i$ are on the
same side of the linear hyperplane in $\RR^\na$ spanned be the vectors $\ga_j$
with
$j\in I\setminus\{i\}\,.$ This implies $x_i>0$ in the expansion
$\xi\,=\,\sum_{j\in I} x_j\ga_j$. So $i\not\in I^-$.

If $I\setminus\{i\}$ is contained in a maximal simplex $I^\prime\neq I$, then
$J\subset I^\prime$ and hence $I<I^\prime$. Now look at the two expansions
$\xi\,=x_i\ga_i+\,\sum_{j\in I\cap I^\prime} x_j\ga_j$ and
$\xi\,=\,y_k\ga_k+\,\sum_{j\in I\cap I^\prime} y_j\ga_j$ where
$\{k\}=I^\prime\setminus (I\cap I^\prime)$. Then $y_k<0$ because
$I^{\prime -}\not\subset I$ by the preceding lemma. On the other hand, $x_i$
and $y_k$ have
different signs because $\ga_i$ and $\ga_k$ lie on different sides of the
linear hyperplane spanned by the vectors $\ga_j$ with $j\in I\cap I^\prime$. We
see $x_i>0$ and $i\not\in I^-$.

Conclusion: $I^-\subset J$. \qed

\begin{proposition}\label{generate}
The elements $c_I$  $(\,I\in\gT^\na\,)$ generate $\cR_{\sA,\gT}$ as a
 $\ZZ[D^{-1}]$-module.
\end{proposition}
\textbf{proof:}
$\cR_{\sA,\gT}$ is linearly generated by monomials
in $c_1,\ldots,c_N\,.$  For one $I_0\in\gT^\na$ we have $I_0^-=\emptyset\,,$
hence $c_{I_0}=1\,.$ Therefore we only need to show that for every $j$ and
every $I_1$ the product $c_j\cdot c_{I_1}$ can be written as a linear
combination of $c_I$'s. If $j\in I_1$ one can use the linear relations
(\ref{eq:linrel}) to express
every $c_i$ with $i\in I_1$ as a $\ZZ[D^{-1}]$-linear combination of $c_k$'s
with $k\not\in I_1\,.$ Since this works for $c_j$ in particular,
the problem can be reduced to showing that a monomial of the form
$\prod_{i\in J} c_i$ with $J$ a simplex of the triangulation,
can be written as a linear combination of $c_I$'s.
Given such a $J$ take $I_J\in\gT^\na$ such that $I_J^-\subset J\subset I_J\,;$
see lemma \ref{tussen}. If
$J=I_J^-$, then  $\prod_{i\in J} c_i\,=\,c_{I_{\!J}}$ and we are done. If
$J\neq I_J^-$ take $m\in J\setminus I_J^-$ and use the linear relations
(\ref{eq:linrel}) to rewrite
$c_m$ as a $\ZZ[D^{-1}]$-linear combination of $c_k$'s with $k\not\in I_J\,.$
This leads to an expression for $\prod_{i\in J} c_i$ as a $\ZZ[D^{-1}]$-linear
combination of monomials of the form
$\prod_{i\in K} c_i$ with $K$ a simplex of the triangulation
and $I_J^-\subsetneq K\,.$ Given such a $K$ take $I_K\in\gT^\na$ such that
$I_K^-\subset K\subset I_K\,.$ Then, according to lemma \ref{grot},
$I_J<I_K\,.$ We proceed by induction.\qed

\

Next we follow Danilov's arguments in
\cite{dani} remark 3.8 to prove
\begin{equation}\label{eq:poinc}
\sum_{k\geqslant 0} \left(\dim _{\QQ}
\,\cR_{\sA,\gT}^{(k)}\otimes\QQ\right)\,\tau^k\:=\:
\sum_{m=0}^{\na}\sharp\, (\,\gT^m\,)\,\tau^m(1-\tau)^{\na-m}
\end{equation}
We have added a few references of which \cite{munk}
is most relevant because it deals with a triangulation of a polytope,
while \cite{dani} deals with a triangulation of a sphere.
In \cite{stan,munk} the Stanley-Reisner ring $\QQ[\gT\,]$
of the simplicial complex $\gT$ over the field $\QQ$
is defined as the quotient of the polynomial ring $\QQ[C_1,\ldots,C_N]$ modulo
the ideal generated by the monomials (\ref{eq:srrel}).
$\QQ[\gT]$ is a Cohen-Macaulay ring of Krull dimension $\na$; see \cite{munk}
thm.2.2 and \cite{stan} thm 1.3.
By definition \ref{ring} there is a natural homomorphism
$\QQ[\gT]\rightarrow\cR_{\sA,\gT}\otimes\QQ$ with kernel generated by
the $\na$ elements $\alpha_i:=a_{i1}C_1+a_{i2}C_2+\ldots+a_{iN}C_N$.
By proposition \ref{generate} the ring $\cR_{\sA,\gT}\otimes\QQ$
is a finite dimensional $\QQ$-vector space and hence has Krull dimension $0$.
It also follows from proposition \ref{generate} that
$\cR_{\sA,\gT}\otimes\QQ$ and $\QQ[\gT]$ are local rings. We can now apply
\cite{matsumura} thm.16.B and see that
$\alpha_1,\ldots,\alpha_\na$ is a regular sequence. As pointed out in
\cite{dani} remark 3.8b this implies that the Poincar\'e series of
$\cR_{\sA,\gT}\otimes\QQ$ is equal to
$(1-\lambda)^\na$ times the Poincar\'e series of $\QQ[\gT\,]$.
The latter is
$\sum_{m=0}^\na\,\sharp\,(\gT^m)\,\lambda^m(1-\lambda)^{-m}$
by \cite{stan} thm. 1.4 (where it is called Hilbert series).
Formula (\ref{eq:poinc}) follows.

We see that $\dim_\QQ \cR_{\sA,\gT}\otimes\QQ=\sharp\,\gT^\na$ and hence that
the elements $c_I$  $(I\in\gT^\na)$ are linearly
independent over $\QQ$.
\emph{This completes the proof of theorem \ref{dims}}\qed

\

\begin{corollary}\label{allvert}
If $\gT$ is unimodular, then
\begin{enumerate}
\item
 $\cR_{\sA,\gT}$ is a free $\ZZ$-module with rank equal to $\vol\Delta$.
\item
$\Delta\,\cap\,\ZZ^\na \:=\:\gT^1\:=\: \{\ga_1,\ldots,\ga_N\}$
\item there is an isomorphism
$\cR_{\sA,\gT}^{(1)}\stackrel{\sim}{\rightarrow} \LL_\ZZ^\vee $
such that $ c_j\leftrightarrow\gb_j$ ($j=1,\ldots,N$)
\end{enumerate}
\end{corollary}
\textbf{proof:} \textbf{(i)} immediately follows from theorem \ref{dims}.
\\
\textbf{(ii)}
Assume that there is a lattice point in $\Delta$ which is not a vertex of
$\gT\,.$ This point lies in
some maximal simplex and gives rise to a decomposition of this simplex into at
least two integral simplices. This contradicts the assumption.
\\
\textbf{(iii)} Because of (ii) all monomials in (\ref{eq:srrel}) have degree
$\geq 2$. Consequently,
$\cR_{\sA,\gT}^{(1)}$ is just the quotient of
$\ZZ\,C_1\oplus\ldots\oplus\ZZ\,C_N$ modulo the span of the linear forms in
(\ref{eq:linrel}). This quotient
is $\LL_\ZZ^\vee$.
\qed


\section{A domain of definition for the function $\Psi_{\gT,\beta}\,.$}
\label{solutions}

We first investigate for which $\lambda$'s one possibly has
$Q_\lambda(\ec)\neq 0$ in $\cR_{\sA,\gT}\otimes\QQ\,.$

For $I\in\gT^\na$ let $\sA_I$ denote the $\na\times\na$-submatrix
of $\sA$ with columns $\ga_i\;(i\in I)\,.$ By
$\ep_I$ we denote the $N\times N$-matrix whose entries with row index not in
$I$ are all $0$ and whose $\na\times N$-submatrix
of entries with row index in $I$ is $\sA_I^{-1}\sA\,.$
This $\ep_I$ is an idempotent linear operator on $\RR^N\,.$ Now define:
\begin{equation}\label{eq:project}
\gP_\gT\::=\:\conv\{\ep_I\:|\:I\in\gT^\na\}
\hspace{4mm}\textrm{in}\hspace{3mm}\textrm{Mat}_{N\times N}(\RR)\,.
\end{equation}
The image of the idempotent operator $\id-\ep_I$ is $\LL_\RR\,.$
Therefore all elements of $\id-\gP_\gT=\conv\{\id-\ep_I\:|\:I\in\gT^\na\}$
are idempotent
operators on $\RR^N$ with image $\LL_\RR\,.$ Hence all elements of
$\gP_\gT$ are also idempotent operators on $\RR^N$.

For every $\lambda\in\ZZ^N$ one has the polytope $\gP_\gT(\lambda)$
which is the
convex hull of $\{\ep_I(\lambda)\:|\:I\in\gT^\na\}$
in $\RR^N\,.$ This obviously
depends only on $\lambda\bmod \LL\,.$

\

\begin{lemma}\label{support}
If $\lambda\in\ZZ^N$ is such that $Q_\lambda(\ec)\neq 0$ in
$\cR_{\sA,\gT}\otimes\QQ$, then  $\lambda$ lies in the set
$\gP_\gT(\lambda)\:+\:\,\cC_\gT^\vee$. Here
$\cC_\gT^\vee$ is the dual of
the cone $\cC_{\gT}$ defined in (\ref{eq:ct}):
\begin{equation}\label{eq:duco}
\cC_{\gT}^\vee\::=\:\{\ell\in\LL_\RR\;|\;\omega\cdot\ell\geq 0
\;\textrm{ for all }\;\omega\in\cC_\gT\;\}\,.
\end{equation}

\end{lemma}
\textbf{proof:}
If $Q_\lambda(\ec)\neq 0$ then $\{i\,|\:\lambda_i<0\}$ is
contained in some maximal simplex $I$ of $\gT\,.$ Let
$\ell=(\id-\ep_I)(\lambda)\,.$
Then $\ell=(\ell_1,\ldots,\ell_N)\in\LL_\RR$ and
$\gb_j\cdot\ell=\ell_j=\lambda_j\geq 0$ for all $j\in I^\ast\,.$
This shows $(\id-\ep_I)(\lambda)\in\cC_{I^\ast}^\vee\subset\cC_\gT^\vee\,.$
\qed

\

\begin{lemma}\label{estimate}
The coefficients in the power series expansion
\begin{equation}\label{eq:qx}
\frac{\prod_{\lambda_j<0}\prod_{k=0}^{-\lambda_j-1}(k+x_j)}
{\prod_{\lambda_j>0}\prod_{k=1}^{\lambda_j}(k-x_j)}
\,=\,\sum_{m_1,\ldots,m_N\geq 0} K_{m_1,\ldots,m_N}
x_1^{m_1}\cdot\ldots\cdot x_N^{m_N}
\end{equation}
satisfy
\begin{equation}\label{eq:schat}
0\leq K_{m_1,\ldots,m_N}\,\leq\, N^{\parallel\lambda\parallel}
\cdot2^{\parallel m\parallel+N}\cdot N!\cdot
(\max(1,N-\deg\lambda))!
\end{equation}
with  $\parallel m\parallel:=\sum_{i=1}^N m_i$
and $\parallel\lambda\parallel:=\sum_{i=1}^N |\lambda_i|$
and $\deg\lambda:=\sum_{i=1}^N \lambda_i=\ga_0^\vee\cdot\beta\,.$
\end{lemma}
\textbf{proof:} Clearly $K_{m_1,\ldots,m_N}\geq 0$. Clearly also $2^{-\parallel
m\parallel}K_{m_1,\ldots,m_N}$ is less than
the value of the left hand side at $x_1=\ldots=x_N\,=\,\frac{1}{2}\,.$
Therefore
$$
K_{m_1,\ldots,m_N}\:<\:
2^{\parallel m\parallel +S}\cdot\frac{\prod_{\lambda_j<0}\,(-\lambda_j)!}{
\prod_{\lambda_j>0}\,(\lambda_j-1)!} \:\leq\:
2^{\parallel m\parallel+S}\cdot\frac{P!}{(R-S)!}\cdot S^{R-S}
$$
where $P\,=\,-\sum_{\lambda_i<0}\lambda_i$ and
$R\,=\,\sum_{\lambda_i>0}\lambda_i$ and $S\,=\,\sharp\{i\,|\,\lambda_i>0\}\,.$

If $P\leq R-S$ then $\frac{P!}{(R-S)!}\leq 1$.
If $P>R-S$ then
$\frac{P!}{(R-S)!}\,\leq\, 2^P\cdot (P-R+S)!$.
Combining these estimates one arrives at (\ref{eq:schat}).
\qed

\

The sum of the series $\sum_{\lambda\in\cA^{-1}(\beta)}\;
Q_\lambda(\ec)\ee^{\tp\ez\cdot\lambda}$ in formula (\ref{eq:form2})
should be computed as a
limit for $L\rightarrow\infty$ of partial sums $\Sigma_L$
taking only terms with $\parallel\lambda\parallel\:\leq L\,.$
These sums only involve $\lambda$'s with
$Q_\lambda(\ec)\neq 0\,.$ According to lemma \ref{support} such a $\lambda$ is
of the form $\lambda=\tilde{\lambda}+\ell$ with $\ell\in\cC_\gT^\vee$
and with $\tilde{\lambda}$ contained in a compact polytope which only depends
on $\beta$. Therefore
$\parallel \lambda\parallel\leq\parallel \ell\parallel +$
some constant which only depends on $\beta$.
Since
$$
Q_\lambda(\ec)=(-1)^{\sharp\{i\,|\,\lambda_i<0\}}
\sum_{m_1,\ldots,m_N\geq 0,\,\parallel m\parallel\leq \na}
(-1)^{\parallel m\parallel} K_{m_1,\ldots,m_N}
c_1^{m_1}\cdot\ldots\cdot c_N^{m_N}
$$
lemma \ref{estimate} shows that the coordinates of $Q_\lambda(\ec)$
with respect to a basis of the vector space
$\cR_{\sA,\gT}\otimes\QQ$ are less than
$N^{\parallel\ell\parallel}$ times some constant which only depends
on $\beta\,.$
Thus one sees that the limit of the partial sums exists if
the imaginary part $\Im\,\ez$ of $\ez$ satisfies
\begin{equation}\label{eq:domain1}
\Im\,\ez\cdot\ell\:>\:
\frac{\log N}{2\pi}\parallel\ell\parallel
\hspace{4mm}\textrm{ for all  } \ell\in\cC_\gT^\vee
\end{equation}
Let $p:\RR^{N\vee}\rightarrow\LL_\RR^\vee$ denote
the canonical projection. If $b\in\LL_\RR^\vee$ is any vector
which satisfies
\begin{equation}\label{eq:bdom}
b\cdot\ell\:>\:\frac{\log N}{2\pi}\parallel\ell\parallel
\hspace{4mm}\textrm{ for all  } \ell\in\cC_\gT^\vee\,,
\end{equation}
then $b\in\cC_\gT$ and every $\ez$ with the property
$ p\,(\Im\,\ez)\in b+\cC_\gT$ satisfies (\ref{eq:domain1}).
Let us therefore define
\begin{equation}\label{eq:domain4}
 \cB_\gT\::=\:
\bigcup_{b\;\textrm{ s.t. (\ref{eq:bdom})}}\;
(\:b\:+\:\cC_\gT\:)
\end{equation}
The above discussion proves:

\

\begin{theorem}\label{defdom}
Formula (\ref{eq:form2}):
$$
\Psi_{\gT,\beta}(\ez)\::=\:
\sum_{\lambda\in\cA^{-1}(\beta)}\hspace{2mm}
Q_\lambda(\ec)\ee^{\tp\ez\cdot\lambda}\cdot\ee^{\tp\ez\cdot\ec}
$$
defines a function with values in
$\cR_{\sA,\gT}\,\otimes\CC$ on the domain
\begin{equation}\label{eq:defdom}
\cV_\gT\::=\{\ez\in\CC^{N\vee}\;|\; p\,(\Im\,\ez)\in\cB_\gT\;\}\,.
\end{equation}
\qed
\end{theorem}

\

In order to have a more global geometric picture of where the domain of
definition of the function $\Psi_{\gT,\beta}$ is situated we
give a brief description of \emph{the toric variety associated with
the pointed secondary fan.}

The pointed secondary fan is a complete fan of strongly convex polyhedral cones
which are generated by vectors from the lattice
$\LL_\ZZ^\vee\,.$
By the general theory of toric varieties \cite{ful,oda} this lattice-fan  pair
gives rise to a toric variety. In the case of $\LL_\ZZ ^\vee$ and the pointed
secondary fan the general construction reads as follows.

For each regular triangulation $\gT$ one has the  cone $\cC_\gT$ in the
secondary fan and one considers  the monoid ring $\ZZ[\LL_\gT]$ of the
sub-monoid $\LL_\gT$ of $\LL\,:$
\begin{equation}\label{eq:lt}
\LL_\gT\::=\LL\cap\cC_\gT^\vee\:=\:
\{\:\ell\in\LL\:|\: \omega\cdot\ell\geq 0\hspace{3mm}
\textrm{ for all }\omega\in\cC_\gT\:\}\,.
\end{equation}
The affine schemes $\cU_\gT\::=\:{\spec\:}\ZZ[\LL_\gT]$ for the various
triangulations naturally glue together to form the toric variety
for the pointed secondary fan.

A complex point of $\cU_\gT$ is just a homomorphism from the
additive monoid $\LL_\gT$ to the multiplicative monoid $\CC\,.$
There is a special point in $\cU_\gT$, namely the homomorphism sending
$0\in\LL_\gT$ to
$1$ and all other elements of $\LL_\gT$ to $0$. A disc of radius $r\,,$
$0<r<1\,,$ about this special point consists of homomorphisms
$\LL_\gT\rightarrow\CC$
with image contained in the disc of radius $r$ in $\CC\,.$

A vector $\ez\in\CC^{N\vee}$ defines the homomorphism
$$
\LL\rightarrow\CC^\ast\,,\;\;\ell\mapsto \ee^{2\pi i \ez\cdot\ell}
$$
and hence a point of the toric variety. The point lies in the
disc of radius $r<1$ about the special point corresponding to a
regular triangulation $\gT$ iff $\Im\ez\cdot\ell>-\frac{\log r}{2\pi}$
holds for every $\ell\in\LL_\gT$. It suffices of course to require this
only for a set of generators of $\LL_\gT$.

If $b$ is in $\cC_\gT$ and $K$ is such that $K>b\cdot\ell$ for all $\ell$ from
a set of generators of $\LL_\gT$, then the set
$\{\:\ez\in\CC^{N\vee} \;|\; p\,(\Im\,\ez)\in b+\cC_\gT\:\}$
contains the intersection of the disc of radius $\exp(-2\pi K)$
with the torus $\Hom{\LL}{\CC^\ast}\,.$

\emph{This shows that the domain of definition of the function
$\Psi_{\gT,\beta}$ is situated about the special point associated
with $\gT$ in the toric variety of the pointed secondary fan.}


\section{The special case $\beta=0$}
\label{betanul}

The function $\Psi_{\gT,0}$ is
invariant under the action of $\MM_\CC^\vee$; see (\ref{eq:character}). So it
is
in fact a function on the domain $\LL_\RR^\vee\:+\:\sqrt{-1} \cB_\gT$ in
$\LL_\CC^\vee\,.$
For  $F\in\Hom{\cR_{\sA,\gT}}{\CC}$ we have the $\CC$-valued function
$F\Psi_{\gT,0}$ on $\LL_\RR^\vee\:+\:\sqrt{-1} \cB_\gT\,.$

\

\begin{lemma}\label{locsys}
If $F\Psi_{\gT,0}$ is the $0$-function on
$\LL_\RR^\vee\:+\:\sqrt{-1} \cB_\gT$
 then $F\,=\,0$.
\end{lemma}
\textbf{proof:}
By lemma \ref{support} the series $\Psi_{\gT,0}$ involves only $\lambda$'s
in $\LL\cap\cC_\gT^\vee$ and $\lambda=0$ is really present
with $Q_0(\ec)=1$. Moreover $\cB_\gT$ is contained in the interior of
$\cC_\gT\,.$ Therefore, if
$F\Psi_{\gT,0}\,=\,0\,,$ then the polynomial function
$$
\sum_{m_1,\ldots,m_N\geq 0}
\frac{(\tp)^{m_1+\ldots+m_N}}{m_1!\cdot\ldots\cdot m_N!}\cdot
F(c_1^{m_1}\cdot\ldots\cdot c_N^{m_N})\cdot z_1^{m_1}\cdot\ldots\cdot z_N^{m_N}
$$
is bounded on an unbounded open domain in $\CC^N$. So, this is the zero
polynomial. Therefore $F(c_1^{m_1}\cdot\ldots\cdot c_N^{m_N})=0$ for all
$m_1,\ldots,m_N\geq 0\,.$
\qed

\

\begin{theorem}\label{isosol}
If $\beta=0$ and $\gT$ is unimodular, then there is an isomorphism:
$$
\Hom{\cR_{\sA,\gT}}{\CC}\stackrel{\sim}{\rightarrow}
\textbf{solution space of (\ref{eq:gkzlin})-(\ref{eq:gkzmon})}\;,
\hspace{5mm}F\mapsto F\Phi_{\gT,0}
\;.
$$
\end{theorem}
\textbf{proof:}
Lemma \ref{locsys} shows that the map is injective.
{}From corollary \ref{allvert}   we know
$\dim\Hom{\cR_{\sA,\gT}}{\CC}\,=\,\vol\Delta$.
Since the triangulation $\gT$ is unimodular, the proof of
\cite{stur} prop.13.15
shows that the normality condition
for the correction in  \cite{gkz1a} to \cite{gkz1} thm. 5 is satisfied.
Therefore the number of linearly independent
solutions of the GKZ system (\ref{eq:gkzlin})-(\ref{eq:gkzmon})
at a generic point equals $\vol\Delta$.
\qed


\section{Triangulations with non-empty core}
\label{coresection}

The intersection of all maximal simplices in a regular triangulation $\gT$ of
$\Delta$ is a remarkable structure. We call it the core of $\gT$.
It is a simplex in the triangulation $\gT\,.$ Since we identify simplices with
their index sets, we view $\core\gT$ also as a subset of $\{1,\ldots,N\}\,.$

\

\begin{definition}\label{tcore}
\hspace{5mm}$\displaystyle{\core\gT\::=\:\bigcap_{I\in\gT^\na}\,I}$
\end{definition}

\

\begin{lemma}\label{corebound}
A simplex which does not contain $\core\gT$ lies in the boundary of $\Delta\,.$
\end{lemma}
\textbf{proof:}
It suffices to prove this for simplices of the form $I\setminus \{j\}$ with
$I\in\gT^\na$ and $j\in\core\gT\,.$
Since every maximal simplex contains $j$, $I$ is the only maximal simplex which
contains $I\setminus \{j\}\,.$ Therefore
 $I\setminus \{j\}$ lies in the boundary of $\Delta\,.$
\qed

\

\begin{lemma}\label{lcore}
\hspace{5mm}$\displaystyle{\core\gT\:=\:\{j\:|\:\ell_j\leq 0 \;
\textrm{ for all }\;\ell\in\cC_{\gT}^\vee\:\} } $
\end{lemma}
\textbf{proof:}
$\supset\,:$ assume $j\not\in\core\gT$, say $j\not\in I$ for some
$I\in\gT^\na$. Then there is a relation
$\ga_j-\sum_{i\in I} x_i\ga_i\,=\,0\,;$
whence an $\ell\in\LL$ with $\ell_j>0$ and $\{i\,|\,\ell_i<0\}\subset I$.
As in the proof of lemma \ref{support} this implies
$\ell\in\cC_{\gT}^\vee\,.$

$\subset\,:$ assume $j\in\core\gT\,.$ First consider an $\ell\in\LL_\RR$ such
that
$\{i\,|\,\ell_i<0\}$ is a simplex.
Let $L=\sum_{\ell_i>0} \ell_i=\sum_{\ell_i<0} -\ell_i\,.$
The relation in (\ref{eq:L}) can be rewritten as
\begin{equation}\label{eq:convell}
\sum_{\ell_i>0} \frac{\ell_i}{L}\ga_i
\:=\:
\sum_{\ell_i<0} \frac{-\ell_i}{L}\ga_i
\end{equation}
Suppose $\ell_j>0$.
Then the simplex $\{i\,|\,\ell_i<0\}$ lies in a boundary
face of $\Delta\,.$
Take a linear functional $F$ whose restriction to  $\Delta$ attains its maximum
exactly on this face. Evaluate $F$ on both sides of (\ref{eq:convell}). The
value on the right hand side is $\max F$, but on the left hand side it is
$<\max F$, because
$F(\ga_j)<\max F$ and $\ell_j>0$. Contradiction!
Therefore we conclude: $\ell_j\leq 0$ if $\ell$ is such that
$\{i\,|\,\ell_i<0\}$ is a simplex. From the constructions in section
\ref{secfan} one sees that $\{i\,|\,\ell_i<0\}$ is a simplex if and only if
$\ell\in\cC_{I^\ast}^\vee$ for some $I\in\gT^\na$; note:
$\ell_j=\gb_j\cdot\ell$. Since
$\cC_{\gT}^\vee$ is the Minkowski sum of the cones $\cC_{I^\ast}^\vee$ with
$I\in\gT^\na$ we finally get: $\ell_j\leq 0$ for every
$\ell\,\in\cC_{\gT}^\vee\,.$
\qed

\

\begin{definition}
\begin{equation}\label{eq:ccore}
c_{\core}\::=\:\prod_{i\in\core\gT} c_i
\end{equation}
\end{definition}

\

\begin{corollary} \label{coreideal}
If $\lambda$ is such that
$\cA\lambda\,=\,\sum_{i\in\core\gT} m_i \ga_i$  with all  $m_i<0$
then $\lambda_i<0$ for every $i\in\core\gT$ and hence
$$
Q_\lambda(\ec)\,\in\,c_{\core}\cR_{\sA,\gT}
$$
\end{corollary}
\textbf{proof:}
Let $\mu=(\mu_1,\ldots,\mu_N)$ be defined by $\mu_i=m_i$ for $i\in\core\gT$
and $\mu_i=0$ for $i\not\in\core\gT$. Then
$\gP_\gT(\lambda)=\gP_\gT(\mu)$ in lemma \ref{support}.
{}From the definitions one sees immediately that
$\gP_\gT(\mu)=\{\mu\}\,.$ The result now follows from lemmas \ref{support} and
\ref{lcore}.

\qed

\

\begin{corollary}\label{corid}
If $\core\gT$ is not empty and
$\beta\,=\,\sum_{i\in\core\gT} m_i \ga_i\;$
 with all $\; m_i<0\;$
then the function $\Psi_{\gT,\beta}$ takes values in the ideal
$c_{\core}\cR_{\sA,\gT}\otimes\CC\,.$
\qed
\end{corollary}

\

\begin{theorem}\label{corsol}
If $\core\gT$ is not empty and $\beta\,=\,\sum_{i\in\core\gT} m_i \ga_i\;$
with all $\; m_i<0\;$
then the linear map
$$
\Hom{c_{\core}\cR_{\sA,\gT}}{\CC}\longrightarrow
\textbf{solution space of (\ref{eq:gkzlin})-(\ref{eq:gkzmon})}\;,
\hspace{5mm}F\mapsto F\Phi_{\gT,\beta}
$$
is injective.
\end{theorem}
\textbf{proof:}
{}From lemma \ref{support} and the proof of corollary \ref{coreideal}
one sees
that the series $\Psi_{\gT,\beta}$ involves only $\lambda$'s
in $\mu+(\LL\cap\cC_\gT^\vee)$ and that $\lambda=\mu$
is really present:
$$
Q_\mu(\ec)=c_{\core}\cdot U\hspace{4mm}\textrm{with}\hspace{4mm}
U\::=\:\prod_{i\in\core\gT} \prod_{k=1}^{-m_i-1}(c_i-k)\,.
$$
The rest of the proof is analogous to the proof of lemma \ref{locsys}.
In particular, if
$F\Psi_{\gT,\beta}$ is the $0$-function on $\cV_\gT$, then
$F(c_{\core}\cdot U\cdot c_1^{n_1}\cdot\ldots\cdot c_N^{n_N})=0$ for all
$n_1,\ldots,n_N\geq 0\,.$
The desired result now follows because $U$ is invertible in the ring
$\cR_{\sA,\gT}\otimes\QQ$.
\qed

\

\

\


\begin{center}\textbf{PART II}\end{center}


\section*{Introduction II}

One aspect of the mirror symmetry phenomenon
(cf. \cite{bundel1,bundel2}) is that (generalized)
Calabi-Yau  manifolds seem to come in pairs $(X,Y)$ with the
geometries of $X$ and  $Y$ related in a beautifully intricate way.
On one side of the mirror - usually called \emph{the B-side} - it is the
geometry of complex structure, of periods of a holomorphic differential form,
of variations of Hodge structure. On the other side -
\emph{the A-side} - it is the geometry of symplectic structure,
of algebraic cycles and of enumerative questions about curves
on the manifold.

Batyrev \cite{bat1} showed that behind many examples of the mirror symmetry
phenomenon one can see a simple combinatorial
duality. Batyrev and Borisov gave a generalization of this combinatorial
duality and
formulated a \emph{mirror symmetry conjecture for generalized Calabi-Yau
manifolds} in arbitrary dimension ( \cite{babo} 2.17).
The fundamental combinatorial structure is a \emph{reflexive Gorenstein cone}.

\begin{definition}\label{gorco}
\textup{( \cite{babo} definitions 2.1-2.8.)}
A cone $\Lambda$ in $\RR^\na$ is called a \emph{Gorenstein cone}
if it is generated, i.e.
\begin{equation}\label{eq:lambda}
\Lambda=\RR_{\geq 0}\ga_1\,+\ldots+\,\RR_{\geq 0}\ga_N\,,
\end{equation}
by a finite set $\{\ga_1,\ldots,\ga_N\}\subset\ZZ^\na$ which satisfies
condition \ref{gkzcond}.
It is called a \emph{reflexive Gorenstein cone} if both $\Lambda$
and its dual $\Lambda^\vee$ are Gorenstein cones,
\begin{equation}\label{eq:dula}
\Lambda^\vee\::=\:\{\:\ey\in\RR^{\na\vee}\:|\:
\forall\,\ex\in\Lambda\,:\:\ey\cdot\ex\,\geq\,0\:\}\,,
\end{equation}
i.e. there should also exist
a vector $\ga_0\in\ZZ^\na$ and a set
$\{\ga_1^\vee,\ldots,\ga_{N^\prime}^\vee\}\subset\ZZ^{\na\vee}$
of generators for $\Lambda^\vee$ such that
$\ga_i^\vee\cdot\ga_0=1$ for $i=1,\ldots,N^\prime$.
The vectors $\ga_0^\vee$ and $\ga_0$ are uniquely determined by
$\Lambda$. The integer $\ga_0^\vee\cdot\ga_0$ is called \emph{ the index} of
$\Lambda\,.$
\end{definition}

\

For a reflexive Gorenstein cone one has one new datum in addition to the data
for GKZ systems; namely $\ga_0$. It has the very important property
\begin{equation}\label{eq:alint}
\textrm{interior}(\Lambda)\cap\ZZ^\na\,=\,\ga_0+\Lambda\,.
\end{equation}

\

\emph{Our aim is to show that in the case of a mirror pair $(X,Y)$
associated with a reflexive Gorenstein cone $\Lambda$ and
a unimodular regular triangulation $\gT$ whose core is not empty
and is not contained in the boundary of $\Delta\,,$
the periods of a holomorphic differential form
on $X$ are given by the function $\Phi_{\gT,-\ga_0}$
which takes values in the ring
$\modquot{\cR_{\sA,\gT}}{\Ann c_{\core}}\otimes\CC$
and that the ring $\modquot{\cR_{\sA,\gT}}{\Ann c_{\core}}$
is isomorphic with a subring of the Chow ring of $Y$.}

\

This project naturally has a B-side and an A-side which we develop
separately in Part II B and Part II A.
Our method puts some natural restrictions on the generality.
For Part II B we must eventually assume that there is a
unimodular triangulation
$\gT$ of the polytope
\begin{equation}\label{eq:delta}
\Delta\::=\:\conv\{\ga_1,\ldots,\ga_N\}=
\{\ex\in\Lambda\:|\:\ga_0^\vee\cdot\ex\,=\,1\:\}\,.
\end{equation}
This restriction which comes from the use of theorem \ref{isosol},
also implies
\begin{equation}\label{eq:a0z}
\ga_0\in\ZZ_{\geq 0}\ga_1+\ldots+\ZZ_{\geq 0}\ga_N\,.
\end{equation}
So, $\Phi_{\gT,-\ga_0}$ is defined in Part I.
For Part II A we must additionally assume that the core
of $\gT$ is not empty and is not
contained in the boundary of $\Delta\,.$


\begin{center}\textbf{PART II B}\end{center}

\section*{Introduction II B}

For a Gorenstein cone $\Lambda$
we denote the monoid algebra $\CC[\Lambda\cap\ZZ^\na]$ by $\cS_\Lambda$
and view it as a subalgebra of the algebra
$\CC[u_1^{\pm 1},\ldots,u_\na^{\pm 1}]$ by identifying
$\sfm=(m_1,\ldots,m_\na)^t\in\Lambda\cap\ZZ^\na$ with the Laurent
monomial $\eu^\sfm\,:=$ $u_1^{m_1}\cdot\ldots\cdot u_\na^{m_\na}$.
For $\sfm\in\ZZ^\na$ we put
$\deg\eu^\sfm\::=\,\deg\sfm\::=\,\ga_0^\vee\cdot\sfm\,.$ Thus
$\cS_\Lambda$ becomes a graded ring.
The scheme
$\PP_\Lambda\::=\:\textrm{Proj }\cS_\Lambda$ is a projective toric variety.
If $\Lambda$ is a reflexive Gorenstein cone, the zero set in $\PP_\Lambda$
of a global section of $\cO_{\PP_\Lambda}(1)$ is called a \emph{generalized
Calabi-Yau manifold} of dimension $\na-2$ ( \cite{babo} 2.15).

The toric variety $\PP_\Lambda$ is a compactification of the
$\na-1$-dimensional torus
\begin{equation}\label{eq:torus2}
\TT\::=\:\modquot{\TTT}{(\ZZ\ga_0^\vee\otimes\CC^{\,\ast})}
\end{equation}
where
\begin{equation}\label{eq:torus1}
\TTT\::=\:\Hom{\ZZ^\na}{\CC^{\,\ast}}
\:=\:\ZZ^{\na\vee}\otimes\CC^{\,\ast}
\end{equation}
is the $\na$-dimensional torus of $\CC$-points
of $\textrm{Spec }\CC[u_1^{\pm 1},\ldots,u_\na^{\pm 1}]\,.$
A global section of $\cO_{\PP_\Lambda}(1)$ is given by a Laurent
polynomial
\begin{equation}\label{eq:laurent}
\es=\sum_{\sfm\in\Delta\cap\ZZ^\na} v_\sfm \eu^\sfm\,.
\end{equation}
with $\Delta$ as in (\ref{eq:delta}). As in  \cite{bat2} we assume from now on

\

\begin{condition}\label{aisall}
\hspace{10mm}$\{\ga_1,\ldots,\ga_N\}=\Delta\cap\ZZ^\na$
\end{condition}

\

The Laurent polynomial $\es$ gives a function on $\TTT$ which
is homogeneous of degree $1$ for the action of
$\ZZ\ga_0^\vee\otimes\CC^{\,\ast}\,.$
Let
\begin{equation}\label{eq:zs}
\sZ_\es\::=\{\textrm{ zero locus of } \;\es\:\}\subset\TT
\end{equation}
Over the complementary set $\TT\setminus\sZ_\es$ there is a section of
$\TTT\rightarrow\TT$ which identifies  $\TT\setminus\sZ_\es$ with the
zero set $\tsZ_{\es-1}$ of $\es-1$ in $\TTT$:
\begin{equation}\label{eq:varieties}
\TT\setminus\sZ_\es\:\simeq\:\tsZ_{\es-1}\subset\TTT
\end{equation}
One may say that according to Batyrev \cite{bat2}
\emph{the geometry on the B-side of mirror symmetry is encoded in the weight
$\na$ part $\cW_\na H^{\na-1}(\TT\setminus\sZ_\es)$ of the Variation of Mixed
Hodge Structure of $H^{\na-1}(\TT\setminus\sZ_\es)$}; the variation comes from
varying the coefficients $v_\sfm$ in (\ref{eq:laurent}).

\begin{remark}\textup{
One usually formulates Mirror Symmetry with on the B-side the Variation
of Hodge Structure on the $d$-th cohomology of a $d$-dimensional
Calabi-Yau manifold. For a CY hypersurface in a toric variety
the Poincar\'e residue mapping gives an isomorphism with the $d+1$-st
cohomology of the hypersurface complement, at least on the primitive parts
(see {} \cite{bat2} prop.5.3). For a CY complete intersection
of codimension $>1$ in a toric variety one needs besides the
Poincar\'e residue mapping also
corollary 3.4 and remark 3.5 in  \cite{babo} to relate the CYCI's
cohomology to the cohomology of the complement of a generalized
Calabi-Yau hypersurface in a toric variety, i.e. to the situation we are
studying in this paper. Our investigations do however also allow on this B-side
of the mirror
generalized Calabi-Yau hypersurfaces which are not related to CY complete
intersections, although on the other A-side we do eventually want
a Calabi-Yau complete intersection (see  \cite{babo} \S 5 for an example
of mirror symmetry with such an asymmetry between the two sides).
}
\end{remark}

\

In \cite{bat2} Batyrev described the weight and Hodge filtrations of this
Variation of Mixed Hodge Structure (VMHS) in terms
of the combinatorics of $\Lambda\,.$ In particular,
$\cW_\na H^{\na-1}(\TT\setminus\sZ_\es)$ corresponds with the ideal
$\CC[\textrm{interior}(\Lambda)\cap\ZZ^\na]$ in
$\cS_\Lambda\,.$ If $\Lambda$ is a reflexive Gorenstein cone of index
$\kappa\,,$ this ideal is the principal ideal generated by $\eu^{\ga_0}$
(cf. (\ref{eq:alint})) and the part of weight $\na$ and Hodge type
$(\na-\kappa,\kappa)$ has dimension $1\,.$

Batyrev \cite{bat2} also showed that the periods of the rational $(\na-1)$-form
\begin{equation}\label{eq:ratform}
\omega_\mu\::=\:\frac{\eu^\mu}{\es^{\deg\mu}}\,
\frac{dt_2}{t_2}\wedge\ldots\wedge\frac{dt_{\na}}{t_{\na}}
\end{equation}
($\mu\in\Lambda\cap\ZZ^\na\,,\; t_2,\ldots,t_{\na}$
coordinates on $\TT$ ) as functions of the coefficients $v_\sfm$ satisfy a GKZ
system of differential
equations (\ref{eq:gkzlin})-(\ref{eq:gkzmon}) with parameters
$\{\ga_1,\ldots,\ga_N\}$ and $\beta=-\mu\,.$
However, \emph{ not all solutions of this system are $\CC$-linear combinations
of the periods of $\omega_\mu\,.$ Theorem \ref{mainthm}
shows precisely which solutions of this system are $\CC$-linear combinations of
the periods of $\omega_\mu$ in case
$\mu\in\ZZ_{\geq 0}\ga_1+\ldots+\ZZ_{\geq 0}\ga_N\,.$ }

The key point of our method is to study the VMHS on
$H^\na(\TTT\rel\tsZ_{\es-1})\,.$
This has the advantage that \emph{if $\ga_1,\ldots,\ga_N$ generate $\ZZ^\na$,
then
$H^\na(\TTT\rel\tsZ_{\es-1})$ is a hypergeometric $\cD$-module} as in
\cite{gkz1} with parameters
$\{\ga_1,\ldots,\ga_N\}$ and $\beta=0\,;$ see
theorem \ref{hygeoDmod}.

\

If $\es$ is $\Lambda$-regular (cf. definition \ref{regular}) there is an exact
sequence of mixed Hodge structures
\begin{equation}\label{eq:exseq}
0\rightarrow H^{\na-1}(\TTT)\rightarrow H^{\na-1}(\tsZ_{\es-1})\rightarrow
H^\na(\TTT\rel\tsZ_{\es-1})\rightarrow H^{\na}(\TTT)\rightarrow 0
\end{equation}
The left hand $0$ results from a theorem of Bernstein-Danilov-Khovanskii
\cite{dk,bat2}. On the right we used $ H^{\na}(\tsZ_{\es-1})=0$
because $\tsZ_{\es-1}$ is an affine variety of dimension $\na-1\,.$ Writing as
usual $\QQ(m)$ for the $1$-dimensional $\QQ$-Hodge structure which is purely of
weight $-2m$ and Hodge type $(-m,-m)$
one has
\begin{equation}\label{eq:hodgetate}
H^{\na-1}(\TTT)\simeq\QQ^\na\otimes\QQ(1-\na)\:,\hspace{5mm}
H^{\na}(\TTT)\simeq\QQ(-\na)\,.
\end{equation}
Morphisms of mixed Hodge structures are strictly compatible with the weight
filtrations ( \cite{del} thm. 2.3.5). Thus the sequence
(\ref{eq:exseq}) in combination with (\ref{eq:varieties})
gives the isomorphisms
\begin{equation}\label{eq:isos}
\cW_i H^{\na-1}(\TT\setminus\sZ_\es)\stackrel{\simeq}{\longrightarrow}
\cW_i H^{\na-1}(\tsZ_{\es-1})\stackrel{\simeq}{\longrightarrow}
\cW_i H^\na(\TTT\rel\tsZ_{\es-1})
\end{equation}
for $ i \leq 2\na-3\,.$
In particular if $\na\geq 3\,,$ the weight $\na$ part relevant for the geometry
on the B-side of mirror symmetry will get a complete and simple description by
our analysis of the GKZ hypergeometric  $\cD$-module
$H^\na(\TTT\rel\tsZ_{\es-1})\,.$

\

\begin{remark}\textup{
Though it plays no role in this paper I want to point out that there is an
interesting relation with recent work of Deninger \cite{den}.
The group $G$ of diagonal $\na\times\na$-matrices with entries
$\pm 1$ acts naturally on $\TTT=\Hom{\ZZ^\na}{\CC^{\,\ast}}\,.$
{}From the inclusion $\imath: \tsZ_{\es-1} \hookrightarrow\TTT$ one gets the
$G$-equivariant map $G\times\tsZ_{\es-1} \rightarrow\TTT\;,$ $\; (g,z)\mapsto
g\cdot\imath (z)$.
Corresponding to this map there is an exact sequence of mixed Hodge structures
with $G$-action analogous to (\ref{eq:exseq}). Taking isotypical parts for the
character
$\det: G\rightarrow \{\pm 1\}$ and using $H^{\na-1}(\TTT)(\det)=0\:,$
$H^{\na}(\TTT)(\det)\stackrel{\simeq}{\rightarrow}H^{\na}(\TTT)$ and
$H^{\na-1}(G\times\tsZ_{\es-1})(\det)\stackrel{\simeq}{\rightarrow}
H^{\na-1}(\tsZ_{\es-1})$ one finds the short exact sequence
\begin{equation}\label{eq:ses1}
0\rightarrow H^{\na-1}(\tsZ_{\es-1})\rightarrow
H^\na(\TTT\rel(G\times\tsZ_{\es-1}))(\det)\rightarrow H^{\na}(\TTT)\rightarrow
0
\end{equation}
see \cite{den} (12). In \cite{den} remark 2.4 Deninger sketches how the
extension (\ref{eq:ses1}) comes from a Steinberg symbol in the group
$K_\na(\tsZ_{\es-1})$ in the algebraic $K$-theory of $\tsZ_{\es-1}$; in our
coordinates (see remark \ref{homocoord3}) this Steinberg symbol reads
\begin{equation}\label{eq:symbol}
\{u_1,u_2,\ldots,u_\na\}\:\in\:K_\na\left(
\modquot{\CC[u_1^{\pm 1},\ldots,u_\na^{\pm 1}]}{(\es-1)}\right)
\end{equation}
The exact sequence (\ref{eq:exseq}) decomposes into two short exact sequences
\begin{eqnarray}
\label{eq:ses2}
0\rightarrow H^{\na-1}(\TTT)\rightarrow H^{\na-1}(\tsZ_{\es-1})
\rightarrow PH^{\na-1}(\tsZ_{\es-1})\rightarrow 0
\\
\label{eq:ses3}
0\rightarrow PH^{\na-1}(\tsZ_{\es-1})\rightarrow
H^\na(\TTT\rel\tsZ_{\es-1})\rightarrow H^{\na}(\TTT)\rightarrow 0
\end{eqnarray}
which define the \emph{primitive part} of cohomology ( \cite{bat2} def. 3.13).
The relation between the various cohomology groups is best displayed in the
following commutative diagram with injective horizontal and surjective vertical
arrows:
\begin{equation}\label{eq:diagram}
\begin{array}{ccccc}
 H^{\na-1}(\TTT)&\rightarrow&H^{\na-1}(\tsZ_{\es-1})&
\rightarrow&H^\na(\TTT\rel(G\times\tsZ_{\es-1}))(\det)\\
&&\downarrow&&\downarrow\\
&&PH^{\na-1}(\tsZ_{\es-1})&\rightarrow&H^\na(\TTT\rel\tsZ_{\es-1})\\
&&&&\downarrow\\&&&&H^{\na}(\TTT)
\end{array}
\end{equation}
With varying coefficients $v_\sfm$ the story plays in the category of
Variations
of Mixed Hodge Structures. With coefficients $v_\sfm$ fixed in some number
field
the story plays in a category of Mixed Motives. A challenge for further
research is to combine these stories and our results on hypergeometric
systems.}
\end{remark}


\section{VMHS associated with a Gorenstein cone}
\label{vmhs}

In this section we prove theorem \ref{hygeoDmod}. This result is
essentially implicitly contained in    \cite{bat2}. Our proof
is mainly a review of constructions and results in \cite{bat2}.

Shifting emphasis from the polytope $\Delta$ to the cone $\Lambda$ we write
$\cS_\Lambda$ (instead of $S_\Delta$ as in \cite{bat2}) for the monoid algebra
$\CC[\Lambda\cap\ZZ^\na]$ viewed as a subalgebra of
$\CC[u_1^{\pm 1},\ldots,u_\na^{\pm 1}]\,.$  The grading is given by
$\deg\eu^\sfm\,=\,\ga_0^\vee\cdot\sfm$ for $\sfm\in\ZZ^\na$.
A homogeneous element $\es$ of degree $1$ in $\cS_\Lambda$ is a Laurent
polynomial as in (\ref{eq:laurent}):
\begin{equation}\label{eq:laurent2}
\es=\sum_{i=1}^N\;v_i\,\eu^{\ga_i}
\end{equation}
with coefficients $v_i\in\CC\,.$
Let $\TTT$, $\TT$, $\sZ_\es$ and $\tsZ_{\es-1}$ be as in
(\ref{eq:torus2})-(\ref{eq:varieties}).

\

\begin{remark}\label{homocoord1}\textup{
When comparing with \cite{bat2} one should keep in mind that in op.cit. $\na$
is the dimension of the polytope $\Delta$ whereas  here $\na$ is the dimension
of the cone $\Lambda$ and the polytope $\Delta$ has dimension $\na-1\,.$ Also
one has to make
the following change of coordinates on $\ZZ^\na$ and $\ZZ^{\na\vee}$.
The  idempotent $\na\times\na$-matrix
$\ga_1\cdot\ga_0^\vee$ gives rise to a direct sum decomposition
$\ZZ^{\na\vee}=\ZZ\ga_0^\vee\oplus \,\Xi$ and thus to a basis
$\{\ga_0^\vee,\,\alpha_2,\ldots,\,\alpha_{\na}\}$ for
$\ZZ^{\na\vee}$. The coordinate change on $\ZZ^\na$ amounts to multiplying
vectors in $\ZZ^\na$ by the matrix $M=(m_{ij})$ with rows
$\ga_0^\vee,\,\alpha_2,\ldots,\,\alpha_{\na}\,.$
In particular, in the new coordinates $\ga_1,\ldots,\ga_N$ all have first
coordinate $1\,.$}

\textup{
The above coordinate change also induces a change of coordinates
on $\TTT\,:$
$u_j=\prod_{i=1}^\na t_i^{m_{ij}}$. The map $\TTT\rightarrow\TT$ is then just
omitting the coordinate $t_1\,.$
In $t$-coordinates $\es$ takes the form $t_1\cdot f$ where $f$ is a Laurent
polynomial in the variables $t_2,\ldots,t_\na\,.$ Thus $\es$ corresponds
with $F_0$ and $\es-1$ with $F$ in \cite{bat2} def. 4.1.
}
\end{remark}
\begin{remark} \label{homocoord3}\textup{
When comparing with \cite{den} one sees again a shift of dimensions
from $\na$ in op. cit. to $\na-1$ here; $T^\na$ with coordinates
$t_1,\ldots,t_\na$ in op. cit. is our $\TT$ with coordinates
$t_2,\ldots,t_\na\,.$
The polynomial $P$ of op. cit. and our $\es$ are related by
$\es=t_1\cdot P$. The identification of $\TT\setminus\sZ_\es$ with
$\tsZ_{\es-1}$ now gives for the Steinberg symbols
$\{P,t_2,\ldots,t_\na\}\,=\,-\{t_1,t_2,\ldots,t_\na\}
\,=\,\{u_1,u_2,\ldots,u_\na\}$ if the coordinates are ordered such that
$\det M=-1\,.$
}
\end{remark}

\

Before we can state Batyrev's results we need some definitions/notations.
\\
\cite{bat2} def. 2.8 defines an ascending sequence of homogeneous ideals
in $\cS_\Lambda$:
\begin{equation}\label{eq:idealweight}
I_\Delta^{(0)}\subset I_\Delta^{(1)}\subset\ldots\subset
I_\Delta^{(\na)}\subset I_\Delta^{(\na+1)}
\end{equation}
where $I_\Delta^{(k)}$ is
generated by the elements $\eu^\sfm$ with $\sfm$ in $\Lambda\cap\ZZ^\na$
but not in any codimension $k$ face of $\Lambda\,;$
in particular
\begin{equation}\label{eq:weightfilt}
I_\Delta^{(0)}=0\:,\hspace{3mm}
I_\Delta^{(1)}=\CC[\mathrm{interior}(\Lambda)\cap\ZZ^\na]\:,\hspace{3mm}
I_\Delta^{(n)}=\cS_\Lambda^+\:,\hspace{3mm}
I_\Delta^{(n+1)}=\cS_\Lambda
\end{equation}
$\cS_\Lambda^+$ is the ideal  in $\cS_\Lambda$ generated
by the monomials of degree $>0\,.$
\\
\cite{bat2} p.379 defines a descending sequence of $\CC$-vector spaces
in $\cS_\Lambda$:
\begin{equation}\label{eq:hodgefilt}
\ldots\supset\cE^{-k}\supset\cE^{-k+1}\supset\ldots
\supset\cE^{-1}\supset\cE^0\supset\cE^1=0
\end{equation}
where $\cE^{-k}$ is spanned by the monomials
$\eu^\sfm$ with $\deg\eu^\sfm\leq k\,.$
\\
\cite{bat2} def. 7.2 defines the differential operators
\begin{equation}\label{eq:difops}
D_i\::=\:u_i\frac{\partial}{\partial u_i}\,+\,
u_i\frac{\partial \es}{\partial u_i}\;,\hspace{5mm} (i=1,\ldots,\na)
\end{equation}
These operate on $\CC[u_1^{\pm 1},\ldots,u_\na^{\pm 1}]\,,$ preserving
$\cS_\Lambda$ and $\cS_\Lambda^+\,.$
\\
\cite{bat2} thm. 4.8 can be used as a definition:
\begin{definition}\label{regular}
$\es$ is said to be \emph{$\Lambda$-regular} if
$\displaystyle{u_1\frac{\partial \es}{\partial u_1},\:
u_2\frac{\partial \es}{\partial u_2}\ldots,\:
u_\na\frac{\partial \es}{\partial u_\na}}$ is a \emph{regular sequence} in
$\cS_\Lambda\,.$
\end{definition}

\begin{theorem}\label{TZ}
\textup{(summary of results in \cite{bat2})}\\
If $\es$ is $\Lambda$-regular, then there is a commutative diagram
\begin{equation}
\label{eq:HTZ}
\begin{array}{ccccc}
\modquot{\cS_\Lambda^+}{\sum_{i=1}^\na D_i \cS_\Lambda^+}&
\stackrel{\simeq}{\rightarrow} & H^{\na-1}(\tsZ_{\es-1}) &
\stackrel{\simeq}{\rightarrow} & H^{\na-1}(\TT\setminus\sZ_\es)\\[.5em]
\downarrow&&\downarrow&&\\[.5em]
\modquot{\cS_\Lambda}{\sum_{i=1}^\na D_i \cS_\Lambda}&
\stackrel{\simeq}{\rightarrow} & H^{\na}(\TTT\rel\tsZ_{\es-1})&&
\end{array}
\end{equation}
in which the horizontal arrows are isomorphisms.
These isomorphisms restrict to the following isomorphisms relating
(\ref{eq:weightfilt}) and (\ref{eq:hodgefilt}) with the weight and Hodge
filtrations on $ H^{\na-1}(\TT\setminus\sZ_\es)$ and
$ H^{\na}(\TTT\rel\tsZ_{\es-1}) \,.$\\
For $k=-1,0,1,\ldots,\na,\na+1\,:$
\[
\begin{array}{lcl}
\textrm{image } I_\Delta^{(k)} \textrm{ in }
\modquot{\cS_\Lambda^+}{\sum_{i=1}^\na D_i \cS_\Lambda^+}
&\stackrel{\simeq}{\rightarrow}&
\cW_{k+\na-1} H^{\na-1}(\TT\setminus\sZ_\es)
\\[.5em]
\textrm{image } \cE^{-k}\cap \cS_\Lambda^+ \textrm{ in }
\modquot{\cS_\Lambda^+}{\sum_{i=1}^\na D_i \cS_\Lambda^+}
&\stackrel{\simeq}{\rightarrow}&
\cF^{\na-k}H^{\na-1}(\TT\setminus\sZ_\es)
\\[.5em]
\textrm{image } I_\Delta^{(k)} \textrm{ in }
\modquot{\cS_\Lambda}{\sum_{i=1}^\na D_i \cS_\Lambda}
&\stackrel{\simeq}{\rightarrow}&
\cW_{k+\na-1} H^{\na}(\TTT\rel\tsZ_{\es-1})
\\[.5em]
\textrm{image } \cE^{-k}\textrm{ in }
\modquot{\cS_\Lambda}{\sum_{i=1}^\na D_i \cS_\Lambda}
&\stackrel{\simeq}{\rightarrow}&
\cF^{\na-k} H^{\na}(\TTT\rel\tsZ_{\es-1})
\end{array}
\]
\end{theorem}
\textbf{proof:}
The statements for $H^{\na-1}(\TT\setminus\sZ_\es)$ are
theorems 7.13, 8.1 and 8.2 in \cite{bat2}. The statements about
$ H^{\na}(\TTT\rel\tsZ_{\es-1})$ can
also be derived with the methods of op. cit., as follows.
Recall that $H^\ast(\TTT\rel\tsZ_{\es-1})$ is the cohomology of the cone of the
natural map of DeRham complexes $\Omega_{\TTT}^\bullet\rightarrow
\Omega_{\tsZ_{\es-1}}^\bullet$ and that this cone complex is in degrees
$i$ and $i+1$
\begin{equation}\label{eq:relcomplex}
\begin{array}{rcccl}
\ldots\rightarrow &
\Omega_{\TTT}^i\oplus\Omega_{\tsZ_{\es-1}}^{i-1}
&\longrightarrow&
\Omega_{\TTT}^{i+1}\oplus\Omega_{\tsZ_{\es-1}}^i
&\rightarrow\ldots\\[.7em]
&(\omega_1,\omega_2)&\mapsto&
(-d\omega_1,d\omega_2+\omega_1|_{\tsZ_{\es-1}})&
\end{array}
\end{equation}
A basis for the
$\CC[u_1^{\pm 1},\ldots,u_\na^{\pm 1}]$-module
$\Omega_{\TTT}^\bullet$ is given by the forms
$\frac{du_{i_1}}{u_{i_1}}\wedge\ldots\wedge\frac{du_{i_r}}{u_{i_r}}$.
Let $\Omega_{\TTT,0}^\bullet$ denote the subgroup of $\Omega_{\TTT}^\bullet$
consisting of the linear combinations of the basic forms
with coefficients in $\CC\,.$ The standard differential $d$ on
$\Omega_{\TTT}^\bullet$ is $0$ on $\Omega_{\TTT,0}^\bullet$.
The inclusion of complexes $\Omega_{\TTT,0}^\bullet\hookrightarrow
\Omega_{\TTT}^\bullet$ is a quasi-isomorphism. So in (\ref{eq:relcomplex})
we may replace $\Omega_{\TTT}^\bullet$ by $\Omega_{\TTT,0}^\bullet$.

For the proof of {} \cite{bat2} thm.7.13 Batyrev uses the $\CC$-linear map
$\cR:\,\cS_\Lambda^+\rightarrow\Omega_{\tsZ_{\es-1}}^{\na-1}\:,$
$\;\cR(\eu^\sfm):=(-1)^{\deg\sfm-1}(\deg\sfm-1)!\,\eu^\sfm
\frac{dt_2}{t_2}\wedge\ldots\wedge\frac{dt_\na}{t_\na}$
(cf. remark \ref{homocoord1} for the $t$-coordinates).
Let us extend this to a $\CC$-linear map
$\cR:\,\cS_\Lambda\rightarrow\Omega_{\TTT,0}^\na
\oplus\Omega_{\tsZ_{\es-1}}^{\na-1}$ by setting $\cR(1)=\left(
\frac{dt_1}{t_1}\wedge\ldots\wedge\frac{dt_\na}{t_\na}\,,\:0\right)\,.$
This induces a surjective linear map
$
\cS_\Lambda\longrightarrow H^\na(\TTT\rel\tsZ_{\es-1})
$
with $\sum_{i=1}^\na D_i \cS_\Lambda^+$ in its kernel.
Note $\displaystyle{D_i(1)=u_i\frac{\partial\es}{\partial u_i}}\,.$
A direct calculation shows for $i=1,\ldots,\na$:
$$
(-1)^{i-1}\cR(t_i\frac{\partial\es}{\partial t_i})
\,=\, d\left(\frac{dt_1}{t_1}\wedge\ldots\wedge
\widehat{\frac{dt_i}{t_i}}\wedge\ldots\wedge
\frac{dt_\na}{t_\na}\,,\:0\right)
$$
in $\Omega_{\TTT,0}^\na
\oplus\Omega_{\tsZ_{\es-1}}^{\na-1}\,.$ Therefore $\cR$ induces a surjective
linear map
$$
\modquot{\cS_\Lambda}{\sum_{i=1}^\na D_i \cS_\Lambda}
\rightarrow H^{\na}(\TTT\rel\tsZ_{\es-1})\,.
$$
A simple dimension count now shows that this is in fact an isomorphism.

The statements about the Hodge filtration and the weight filtration
on $ H^{\na}(\TTT\rel\tsZ_{\es-1})$ follow from the corresponding
statements for  $H^{\na-1}(\TT\setminus\sZ_\es)$ and from
(\ref{eq:hodgetate}).
\qed

\

The \emph{principal $\sA$-determinant} of Gel'fand-Kapranov-Zelevinskii
\cite{gkz2} is a polynomial $E_{\sA} (v_1,\ldots,v_N)\in
\ZZ[v_1,\ldots,v_N]$ such that (see \cite{bat2} prop. 4.16):
\begin{equation}\label{eq:disc}
\es \textrm{ is } \Lambda\textrm{-regular}
\hspace{5mm}\Longleftrightarrow \hspace{5mm}
E_{\sA} (v_1,\ldots,v_N)\neq 0
\end{equation}
\emph{Now we want to vary the coefficients $v_i$
in (\ref{eq:laurent2}) and work over the ring}
\begin{equation}\label{eq:basering}
\CC[\ev]\::=\:\CC[v_1,\ldots,v_N,E_{\sA}^{-1}]\,.
\end{equation}

Let $\Omega^\bullet$ resp. $\tOmega^\bullet$  denote the
DeRham complex of
$\CC[u_1^{\pm 1},\ldots,u_\na^{\pm 1}]\otimes\CC[\ev]$
relative to $\CC[\ev]$ resp. relative to $\CC\,.$
Define on these complexes a new differential
\begin{eqnarray}
\nonumber
&&\delta: \Omega^i\rightarrow\Omega^{i+1}
\textrm{  resp.  }
\tOmega^i\rightarrow \tOmega^{i+1}
\\
\label{eq:newd}
&&\delta\omega\::=\: d\omega+ d\es\wedge\omega
\end{eqnarray}
where $d$ is the ordinary differential on DeRham complexes.

As a basis for the
$\CC[u_1^{\pm 1},\ldots,u_\na^{\pm 1}]\otimes\CC[\ev]$-module
$\Omega^1$ (resp. $\tOmega^1$ ) we take
$\frac{du_1}{u_1},\ldots,\frac{du_\na}{u_\na}$
(resp. $\frac{du_1}{u_1},\ldots,\frac{du_\na}{u_\na},dv_1,\ldots,dv_N$)
and extend it by taking wedge products to a basis
for $\Omega^\bullet$ (resp. $\tOmega^\bullet$ ).
Let $\Omega_\Lambda^\bullet$ (resp. $\Omega_{\Lambda^+}^\bullet$ ) denote the
subgroups of $\Omega^\bullet$ consisting of
the linear combinations of the given basic forms with coefficients in
$\cS_\Lambda\otimes\CC[\ev]$ (resp. $\cS_\Lambda^+\otimes\CC[\ev]$ ).
Define $\tOmega_\Lambda^\bullet$ (resp. $\tOmega_{\Lambda^+}^\bullet$ )
in the same way as subgroups of $\tOmega^\bullet\,.$
The differential $\delta$ (\ref{eq:newd}) preserves these subgroups.
Thus we get the two complexes
\begin{eqnarray*}
(\Omega_\Lambda^\bullet,\delta)&:&\hspace{3mm}
\Omega_\Lambda^0\,\stackrel{\delta}{\rightarrow}
\Omega_\Lambda^1\,\stackrel{\delta}{\rightarrow}
\ldots\stackrel{\delta}{\rightarrow}
\Omega_\Lambda^{\na-1}\,\stackrel{\delta}{\rightarrow}
\Omega_\Lambda^\na
\\
(\tOmega_\Lambda^\bullet,\delta)&:&\hspace{3mm}
\tOmega_\Lambda^0\,\stackrel{\delta}{\rightarrow}
\tOmega_\Lambda^1\,\stackrel{\delta}{\rightarrow}
\ldots\stackrel{\delta}{\rightarrow}
\tOmega_\Lambda^{\na-1}\,\stackrel{\delta}{\rightarrow}
\tOmega_\Lambda^{\na}
\stackrel{\delta}{\rightarrow}
\tOmega_\Lambda^{\na+1}\,\stackrel{\delta}{\rightarrow}
\ldots\stackrel{\delta}{\rightarrow}\tOmega_\Lambda^{N+\na}
\end{eqnarray*}
Then
\begin{equation}
\label{eq:dmod1}
H^\na(\Omega_\Lambda^\bullet,\delta)\:=\:
\left(\modquot{\cS_\Lambda}{\sum_{i=1}^\na D_i \cS_\Lambda}\right)
\otimes\CC[\ev]
\end{equation}
The Gauss-Manin connection
\begin{equation}
\label{eq:gaussmanin}
\nabla\,:\;H^\na(\Omega_\Lambda^\bullet,\delta)\rightarrow
H^\na(\Omega_\Lambda^\bullet,\delta)\otimes\Omega^1_{\CC[\ev]/\,\CC}
\end{equation}
on this module is described by the Katz-Oda construction
(cf. \cite{katz} \S 1.4) as follows. Lift the given
$\xi\in H^\na(\Omega_\Lambda^\bullet,\delta)$ to an element
$\tilde{\xi}$ in $\tOmega_\Lambda^{\na}$. Then $\nabla\xi$ is the cohomology
class of $\delta\tilde{\xi}\in\tOmega_\Lambda^{\na+1}$ in
$H^\na(\Omega_\Lambda^\bullet,\delta)\otimes\Omega^1_{\CC[\ev]/\,\CC}\,.$
Having $\nabla\xi$ one defines
$\frac{\partial}{\partial v_j}\xi\in H^\na(\Omega_\Lambda^\bullet,\delta)$
by
\begin{equation}\label{eq:deriv}
\nabla\xi\:=\:\sum_{j=1}^N\:
\left(\frac{\partial}{\partial v_j}\xi\right)\,\otimes\,dv_j
\end{equation}
In particular for $\mu\in\Lambda\cap\ZZ^\na$ and
\begin{equation}\label{eq:ximu}
\xi_\mu\::=\:\textrm{cohomology class of }
\eu^\mu\cdot\frac{du_1}{u_1}\wedge\ldots\wedge\frac{du_\na}{u_\na}\;
\in H^\na(\Omega_\Lambda^\bullet,\delta)
\end{equation}
we find
\begin{eqnarray}\label{eq:dximu1}
\frac{\partial}{\partial v_j}\xi_\mu&=&\textrm{cohomology class of }
\eu^{\ga_j+\mu}\cdot\frac{du_1}{u_1}\wedge\ldots\wedge\frac{du_\na}{u_\na}
\\  \label{eq:dximu2}
&=& \xi_{\mu+\ga_j}
\end{eqnarray}
The form $\xi_\mu$ for $\mu\neq 0$  corresponds via (\ref{eq:dmod1})
and  {} \cite{bat2} thm.7.13 with the form $\omega_\mu$ in
(\ref{eq:ratform}); more precisely
 $\xi_\mu$ is the cohomology class of  $\omega_\mu$ modulo
$H^{\na-1}(\TTT)$.

\

\begin{corollary}\label{gkzximu}
\begin{eqnarray}
\label{eq:gkzxilin}
\left(\mu+
\sum_{j=1}^N\:\ga_j\,v_j\frac{\partial}{\partial v_j}\right)\;\xi_\mu &=&0
\\
\label{eq:gkzximon}
\left(\prod_{\ell_j>0} \left[\frac{\partial}{\partial v_j}\right]^{\ell_j}
\:-\:\prod_{\ell_j<0} \left[\frac{\partial}{\partial v_j}\right]^{-\ell_j}
\right)\xi_\mu&=&0\hspace{3mm}\textrm{ for } \ell\in\LL
\end{eqnarray}
\end{corollary}
\textbf{proof:}
On the level of differential forms in the complex
$(\Omega_\Lambda^\bullet,\delta)$
the $i$-th equation of (\ref{eq:gkzxilin}) reads
\begin{eqnarray*}
&&\left(\mu_i+
\sum_{j=1}^N\:a_{ij}\,v_j\frac{\partial}{\partial v_j}\right)
\eu^\mu\cdot\frac{du_1}{u_1}\wedge\ldots\wedge\frac{du_\na}{u_\na}\;=
\\
&&=\;
\delta\left((-1)^{i-1}\eu^\mu
\frac{du_1}{u_1}\wedge\ldots\wedge\frac{du_{i-1}}{u_{i-1}}\wedge
\frac{du_{i+1}}{u_{i+1}}\wedge\ldots\wedge\frac{du_\na}{u_\na}\right)
\end{eqnarray*}
(\ref{eq:gkzximon}) follows immediately from (\ref{eq:dximu1}).
\qed

\

\begin{remark}\textup{
We have essentially repeated the proof of   \cite{bat2} thm. 14.2.
There is however a small difference: Batyrev uses coefficients in
$\cS_\Lambda^+$ where we are using coefficients in $\cS_\Lambda\,.$
His differential equations hold for $H^{\na-1}(\TT\setminus\sZ_\es)
=H^n(\Omega_{\Lambda^+}^\bullet,\delta)$
whereas ours only hold in the primitive part
$PH^{\na-1}(\TT\setminus\sZ_\es)\,.$
On the other hand we can also treat $\xi_0$. The following theorem
shows that this gives an important advantage.}
\end{remark}

\begin{theorem}\label{hygeoDmod}
If $\Lambda\cap\ZZ^\na\:=\:\ZZ_{\geq 0}\ga_1+\ldots+\ZZ_{\geq 0}\ga_N\,,$
then  $\xi_0$ generates $H^\na(\Omega_\Lambda^\bullet,\delta)$  as a module
over the ring $\cD\::=\:\CC[v_1,\ldots,v_N,E_{\sA}^{-1},
\frac{\partial}{\partial v_1},\ldots,\frac{\partial}{\partial v_N}]\,.$
The annihilator of $\xi_0$ in \rule{0mm}{3mm}$\cD$ is the left ideal generated
by
the differential operators
\[
\sum_{j=1}^N\:a_{ij}\,v_j\frac{\partial}{\partial v_j}
\hspace{4mm}\textrm{and} \hspace{4mm}
\prod_{\ell_j>0} \left[\frac{\partial}{\partial v_j}\right]^{\ell_j}
\:-\:\prod_{\ell_j<0} \left[\frac{\partial}{\partial v_j}\right]^{-\ell_j}
\]
with $1\leq i\leq\na$ and $\ell\in\LL\,.$

In other words, $ H^{\na}(\TTT\rel\tsZ_{\es-1})=
H^\na(\Omega_\Lambda^\bullet,\delta)$ is the
\emph{hypergeometric $\cD$-module in the sense of {} \cite{gkz1}} \S 2.1 with
parameters
$\{\ga_1,\ldots,\ga_N\}$ and $\beta=0\,.$
\end{theorem}
\textbf{proof:}
Let $\cM_0$ denote the hypergeometric $\cD$-module with parameters
$\beta=0$ and $\{\ga_1,\ldots,\ga_N\}$ as in \cite{gkz1} section 2.1. By
corollary \ref{gkzximu} and formula (\ref{eq:dximu2}) we have a surjective
homomorphism of
$\cD$-modules  $\cM_0\rightarrow H^\na(\Omega_\Lambda^\bullet,\delta)\,.$
The filtration of $\cD$ by the order of differential operators induces an
ascending filtration on $\cM_0$ and $H^\na(\Omega_\Lambda^\bullet,\delta)\,.$
It suffices to prove that the above surjection induces an isomorphism for the
associated graded modules. According to {} \cite{gkz1} prop.3
$gr\,\cM_0$ is isomorphic to the quotient of the ring
$\CC[x_1,\ldots,x_N]\otimes\CC[\ev]$ by the ideal generated by the
linear forms $\sum_{j=1}^N\,a_{ij}x_j$ for $i=1,\ldots,\na$ and by the
polynomials
$\prod_{\ell_j>0} x_j^{\ell_j}\,-\,\prod_{\ell_j<0} x_j^{-\ell_j}$ with
$\ell\in\LL\,.$ Via the substitution homorphism $x_j\mapsto \eu^{\ga_j}$
this quotient ring is isomorphic to the quotient of the ring
$\cS_\Lambda\otimes\CC[\ev]$ by the ideal generated by
$\displaystyle{u_1\frac{\partial \es}{\partial u_1},\:
u_2\frac{\partial \es}{\partial u_2}\,,\ldots,\:
u_\na\frac{\partial \es}{\partial u_\na}}\,.$ Using
(\ref{eq:dximu1}), (\ref{eq:dmod1}) and (\ref{eq:difops}) one checks that the
latter quotient ring is isomorphic to
$gr\,H^\na(\Omega_\Lambda^\bullet,\delta)\,.$
\qed

\

\begin{center}\textbf{PART II A}\end{center}

\section*{Introduction II A}

In this Part II A we give our results the flavor of Mirror Symmetry
by showing that for a regular triangulation $\gT$
which satisfies conditions  (\ref{eq:core1}), (\ref{eq:core2}),
(\ref{eq:vol1}),
 the ring $\cR_{\sA,\gT}$ is the cohomology ring of a
toric variety constructed somehow from the dual Gorenstein cone $\Lambda^\vee$
and that the ring
$ \modquot{\cR_{\sA,\gT}}{\Ann c_{\core}}$
is a subring of the Chow ring of a Calabi-Yau complete intersection
in that toric variety; more precisely the subring is the image of the Chow ring
of the ambient toric variety.

We construct several toric varieties which are also used in \cite{babo}. As we
want to promote the use of triangulations we give a construction of these toric
varieties as a quotient of an open part of
$\CC^d$ ($d$ an appropriate dimension) by a torus. The torus is related to
$\LL$ and the open part is given by the triangulation $\gT$.
Such a construction  of toric varieties is well known (see for instance
\cite{guillemin}).

\section{Triangulations with non-empty core and
 completely split reflexive Gorenstein cones.}
\label{gorenstein}

\begin{proposition}\label{corecs}
Assume that $\gT$ satisfies the following three conditions
\begin{eqnarray}
&&\core\gT\textrm{ is not empty and }
\core\gT\,=\,\{1,\ldots,\kappa\}\label{eq:core1}\\
&&\core\gT\textrm{ is not contained in the boundary of } \Delta
\label{eq:core2}\\
&&\gT \textrm{ is unimodular}\label{eq:vol1}
\end{eqnarray}
Then
$\Lambda\::=\:\RR_{\geq 0}\ga_1\,+\ldots+\,\RR_{\geq 0}\ga_N$
is a reflexive Gorenstein cone of index $\kappa$ and
the dual cone $\Lambda^\vee$ is \emph{completely split} in the
sense of  {} \cite{babo} definition 3.9.
\end{proposition}
\textbf{proof:}
By lemma \ref{corebound} and hypotheses (\ref{eq:core1}) and (\ref{eq:core2})
the $(\na-2)$-dimensional simplices in the boundary of $\Delta$ are precisely
the simplices $I\setminus\{i\}$ with $I\in\gT^\na$ and $i=1,\ldots,\kappa$.
It follows that the dual cone $\Lambda^\vee$ is generated by the set of
row vectors $\{\ga_{I,i}^\vee\;|\;I\in\gT^\na,\;i=1,\ldots,\kappa\:\}$
where
\[
\ga_{I,i}^\vee\;:=\;\textrm{ the } i\textrm{-th row of the matrix }
\sA_I^{-1}
\]
Hypothesis (\ref{eq:vol1}) implies $\ga_{I,i}^\vee\in\ZZ^{\na\vee}$ for all
$I,i\,.$
By construction
\begin{equation}\label{eq:splitting}
\ga_{I,i}^\vee\cdot\ga_j\:=\:
\left\{\begin{array}{ll}\geq  0 & \textrm{ for } j=1,\ldots,N\\
 1 & \textrm{ if } j=i\\
 0 & \textrm{ if } 1\leq j\leq\kappa\,,\;j\neq i
\end{array}\right.
\end{equation}
So if we take
\begin{equation}\label{eq:a0}
\ga_0\::=\:\ga_1+\ldots+\ga_\kappa\;\in\RR^\na
\end{equation}
then
$$
\ga_{I,i}^\vee\cdot\ga_0\,=\,1\hspace{5mm}\textrm{for}\hspace{3mm}
I\in\gT^\na,\;i=1,\ldots,\kappa\,.
$$
This shows that $\Lambda^\vee$ is a Gorenstein cone. Hence $\Lambda$
is a reflexive Gorenstein cone with index $\ga_0^\vee\cdot\ga_0=\kappa\,.$

Every element of $\Lambda^\vee$ can be written as
$\sum_{I,i}\:s_{I,i}\ga_{I,i}^\vee$ with all
$s_{I,i}\in\RR_{\geq 0}\,.$ Such a sum can be rearranged as
$\sum_{i=1}^\kappa t_i \alpha_i$
with $t_i=\sum_I\,s_{I,i}$ and $\alpha_i\in\square_i$ where
\begin{equation}\label{eq:boxi}
\square_i\::=\,\conv\{\ga_{I,i}^\vee\:|\:I\in\gT^\na\:\}\,.
\end{equation}
$\square_i$ is a lattice polytope in the $(\na-\kappa)$-dimensional affine
subspace of $\RR^{\na\vee}$ given by the equations $\xi\cdot\ga_i=1$
and $\xi\cdot\ga_j=0$ if $1\leq j\leq\kappa\,,\;j\neq i$
(cf. (\ref{eq:splitting})).
This shows that $\Lambda^\vee$ is a \emph{completely split} reflexive
Gorenstein cone of index $\kappa$ in the sense of {} \cite{babo}
definition 3.9.

Note that the dimension of $\square_i$ equals $\na-2$ minus
the dimension of the minimal
face of $\Delta$ which contains $\{\ga_j\:|\:j\in\core\gT\setminus\{i\}\;\}\,.$
\qed

\section{Triangulations and toric varieties}
\label{tritorvar}

\emph{We assume from now on that $\gT$ satisfies the conditions
(\ref{eq:core1}), (\ref{eq:core2}), (\ref{eq:vol1}).}

Take some $I_0\in\gT^\na$ and
consider the matrix $(u_{ij})\,:=\,\sA _{I_0}^{-1}\sA$. Then in definition
\ref{ring} the linear forms
\begin{equation}\label{eq:linrelu}
u_{i1}C_1+u_{i2}C_2+\ldots+u_{iN}C_N\hspace{5mm}(i=1,\ldots,\na)
\end{equation}
together with the monomials in (\ref{eq:srrel}) give another system of
generators for the ideal $\cJ$. The corresponding relations in $\cR_{\sA,\gT}$
for $i=1,\ldots,\kappa$ express $c_1,\ldots,c_\kappa$ as linear combinations of
 $c_{\kappa+1},\ldots,c_N$. The relations for $i=\kappa+1,\ldots,\na$
do not involve  $c_1,\ldots,c_\kappa\,.$
Also the monomials in (\ref{eq:srrel}) do not involve $C_1,\ldots,C_\kappa\,.$

Let $\gu_{\kappa+1},\ldots,\gu_N\,\in\RR^{\na-\kappa}$ be the columns
of the matrix $(u_{ij})_{\kappa<i\leq\na,\kappa<j\leq N}$.
There is a simplicial fan $\cF^\prime$ in $\RR^{\na-\kappa}$ given by the cones
\begin{equation} \label{eq:profan}
\RR_{\geq 0}\gu_{i_1}+\ldots+\RR_{\geq 0}\gu_{i_s}\hspace{5mm}
\textrm{ with }i_1,\ldots,i_s > \kappa
\textrm{  and  }\{i_1,\ldots,i_s\}\in\gT
\end{equation}
i.e. the index set is a simplex in the triangulation $\gT$.

The fan $\cF^\prime$ is complete iff $0\in\RR^{\na-\kappa}$ is a linear
combination
with positive coefficients of the vectors $\gu_{\kappa+1},\ldots,\gu_N$.
This is equivalent to condition (\ref{eq:core2}).
Condition (\ref{eq:vol1}) implies that $\cF^\prime$ is a fan of regular
simplicial cones, i.e. its maximal cones are spanned by a basis of
$\ZZ^{\na-\kappa}$.

Combining these considerations with \cite{dani} thm.10.8 or \cite{ful}
prop.p.106 we find:

\begin{theorem}\label{protor}
If the triangulation $\gT$ satisfies  conditions (\ref{eq:core1}),
(\ref{eq:core2}), (\ref{eq:vol1}), then
$\cR_{\sA,\gT}$ is isomorphic to the cohomology ring $H^\ast(\PP_\gT,\ZZ)$
of the $(\na-\kappa)$-dimensional smooth projective toric variety $\PP_\gT$
associated with the fan  $\cF^\prime$ (see definition \ref{ptor}); more
precisely:
$$
\cR_{\sA,\gT}^{(m)}\simeq H^{2m}(\PP_\gT,\ZZ)\;,\hspace{5mm}
m=0,1,\ldots,\na-\kappa\,.
$$
and $\cR_{\sA,\gT}^{(m)}=0$ for $m>\na-\kappa\,.$\qed
\end{theorem}

\

There is much more geometry in those three conditions than was used for
theorem \ref{protor}.
Consider in $\RR^\na$ the fan $\cF$ consisting of the cones
\begin{equation} \label{eq:fan}
\RR_{\geq 0}\ga_{i_1}+\ldots+\RR_{\geq 0}\ga_{i_s}\;\;,\hspace{5mm}
\{i_1,\ldots,i_s\}\in\gT\,.
\end{equation}
The standard constructions produce a toric variety $\EE_{\,\gT}$ from this fan.
We recall the construction of the toric variety $\EE_{\,\gT}$ as a quotient
of an open part of $\CC^N$ by the torus $\LL\otimes\CC^{\,\ast}$. This torus
appears here because $\LL$ is the lattice of linear relations between the
vectors $\ga_1,\ldots,\ga_N\,;$ by condition (\ref{eq:vol1}) and corollary
\ref{allvert} these are exactly the generators of the $1$-dim cones of the fan
$\cF\,.$

Take $\CC^N$ with coordinates $x_1,\dots,x_N$ and define
\begin{eqnarray}
\CC^N_I&\::=\:&\{\,(x_1,\dots,x_N)\in\CC^N\;|\;x_j\neq 0\textrm{  if  }
j\not\in I\:\}\hspace{3mm}\textrm{ for }\hspace{3mm}I\in\gT^\na
\nonumber \\
\label{eq:cnt}
\CC^N_\gT&\::=\:&\bigcup_{I\in\gT^\na}\:\CC^N_I
\end{eqnarray}
The torus $\CC^{\,\ast N}$ acts on $\CC^N$ via coordinatewise multiplication.
The inclusion $\LL\subset\ZZ^N$ induces an inclusion of tori
$\LL\otimes\CC^{\,\ast}\subset\CC^{\,\ast N}$. Thus
$\LL\otimes\CC^{\,\ast}$
acts on $\CC^N\,.$
For $\ell=(\ell_1,\ldots,\ell_N)\in\LL\,,\;t\in\CC^{\,\ast}$ the element
$\ell\otimes t$ acts as
\begin{equation}\label{eq:action}
(\ell\otimes t)\cdot(x_1,\dots,x_N)\::=\:
(t^{\ell_1}x_1,\ldots,t^{\ell_N}x_N)
\end{equation}

\begin{definition}\label{etor}\hspace{5mm}
$\displaystyle{\EE_{\,\gT}\::=\:
\modquot{\CC^N_\gT}{\LL\otimes\CC^{\,\ast}}\,.}$
\end{definition}

\

Take an $\bnl\times N$-matrix $\sB$ with entries in $\ZZ$ such that the columns
of $\sB^t$ constitute a basis for $\LL\,.$
For $I\subset\{1,\ldots,N\}$ we denote by $\sA_I$ (resp.
$\sB_{I^\ast}$) the submatrix of $\sA$  (resp. $\sB$) composed of the entries
with column index in $I$ (resp. in $I^\ast\::=\,\{1,\ldots,N\}\setminus I$ ).
Consider  $I=\{i_1,\ldots,i_\na\}\in\gT^\na$. Then $\det(\sB_{I^\ast})\,=\,\pm
\det(\sA_I)\,=\,\pm 1$ by condition (\ref{eq:vol1}).  So $\sB_{I^\ast}$
is invertible over $\ZZ\,.$ From this one easily sees that there
is an isomorphism
\begin{eqnarray}\label{eq:chart}
\CC^\na&\stackrel{\simeq}{\longrightarrow}&
\modquot{\CC^N_I}{\LL\otimes\CC^{\,\ast}}\\
\nonumber
(y_1,\ldots,y_\na)&\mapsto& (x_1,\dots,x_N)
\textrm{ with  } x_j=\left\{\begin{array}{ll} y_t&\textrm{ if } j=i_t\in I\\
1&\textrm{ if } j\not\in I \end{array}\right.
\end{eqnarray}
Hence $\EE_{\,\gT}$ is a smooth toric variety.
The torus
$\modquot{\CC^{\ast N}}{\LL\otimes\CC^{\,\ast}}=\MM\otimes\CC^{\,\ast}$
acts on $\EE_{\,\gT}$ and the variety
$\EE_{\,\gT}$ contains $\MM\otimes\CC^{\,\ast}$ as a dense open subset.

One constructs in the same way the toric variety $\PP_\gT$ from the fan
$\cF^\prime$ (see (\ref{eq:profan})). Now the lattice of linear relations
between the generators $\gu_{\kappa+1},\ldots,\gu_N$ of the $1$-dimensional
cones of the fan $\cF^\prime$ is the image of the composite map
$\LL\hookrightarrow\ZZ^N\twoheadrightarrow \ZZ^{N-\kappa}$. This map
$\LL\rightarrow \ZZ^{N-\kappa}$ is also injective. Take  $\CC^{N-\kappa}$ with
coordinates $x_{\kappa+1},\dots,x_N$ and define
\begin{eqnarray}
\CC^{N-\kappa}_I&\::=\:&\{\,(x_{\kappa+1},\dots,x_N)\in\CC^N\;|\;x_j\neq
0\textrm{  if  } j\not\in I\:\}\hspace{3mm}\textrm{ for
}\hspace{3mm}I\in\gT^\na
\nonumber \\
\CC^{N-\kappa}_\gT&\::=\:&\bigcup_{I\in\gT^\na}\:\CC^{N-\kappa}_I
\end{eqnarray}
$\LL\otimes\CC^{\,\ast}$ is a subtorus of $\CC^{\,\ast {N-\kappa}}$
and acts accordingly; i.e. as in (\ref{eq:action}) using only the coordinates
with index $>\kappa\,.$

\

\begin{definition}\label{ptor}\hspace{5mm}
$\displaystyle{\PP_\gT\::=\:
\modquot{\CC^{N-\kappa}_\gT}{\LL\otimes\CC^{\,\ast}}\,.}$
\end{definition}

\

$\PP_\gT$ is a smooth projective toric variety: smooth for the same reason as
$\EE_{\,\gT}$ and projective because the fan $\cF^\prime$ is complete.
Projection onto the last $N-\kappa$ coordinates induces a surjective morphism
\begin{equation}\label{eq:morphism}
\pi\,:\;\EE_{\,\gT}\rightarrow\PP_\gT
\end{equation}
As (\ref{eq:cnt}) puts no restriction on the coordinates
$x_1,\ldots,x_\kappa\,,$ the fibers of $\pi$ are complex vector spaces
of dimension $\kappa\,;$ more precisely, (\ref{eq:chart}) gives a
trivialization
\[
\modquot{\CC^N_I}{\LL\otimes\CC^{\,\ast}}\;\simeq\;\CC^\na\;\simeq\;
\CC^\kappa\times\CC^{\na-\kappa}\;\simeq\;
\CC^\kappa\times\left(\modquot{\CC^{N-\kappa}_I}{\LL\otimes\CC^{\,\ast}}
\right)
\]
Thus:

\begin{proposition} $\EE_{\,\gT}$ has the structure of a vector bundle of rank
$\kappa$ over $\PP_\gT\,.$
\qed
\end{proposition}

\

The dual vector bundle $\EE_{\,\gT}^\vee\rightarrow\PP_\gT$ can be constructed
as
\begin{equation}\label{eq:edual}
\EE_{\,\gT}^\vee\::=\:
\modquot{\CC^N_\gT}{(\LL\otimes\CC^{\,\ast})^\prime}
\end{equation}
with $\CC^N_\gT$ as in definition \ref{etor}, but with the action of
$\LL\otimes\CC^{\,\ast}$ slightly modified from (\ref{eq:action}):
the element $\ell\otimes t$ now acts as
\begin{equation}\label{eq:duaction}
(\ell\otimes t)\cdot^\prime\,(x_1,\dots,x_N)\::=\:
(t^{-\ell_1}x_1,\ldots,t^{-\ell_\kappa}x_\kappa,
t^{\ell_{\kappa+1}}x_{\kappa+1},\ldots,t^{\ell_N}x_N)
\end{equation}

For the sake of completeness we also describe the construction
of the bundle of projective spaces
$\PP\EE_{\,\gT}\rightarrow\PP_\gT$  associated with the vector bundle
$\EE_{\,\gT}\rightarrow\PP_\gT\,.$
Take as before $\CC^N$ with coordinates $x_1,\dots,x_N\,.$ Define for
$i\in\core\gT$ and $I\in\gT^\na$
\begin{eqnarray}
\CC^N_{i,I}&\::=\:&\{\,(x_1,\dots,x_N)\in\CC^N\;|\;x_i\neq 0 \textrm{ and
}x_j\neq 0\textrm{  if  } j\not\in I\:\}
\nonumber \\
\label{eq:cnpt}
\CC^N_{\gT\circ}&\::=\:&\bigcup_{i\in\core\gT\,,\:I\in\gT^\na}\:\CC^N_{i,I}
\end{eqnarray}
Write
$\ek:=(k_1,\ldots,k_N)^t$ with $k_j=1$ if  $j\in\core\gT$
resp. $k_j=0$ if $j\not\in\core\gT\,,$
i.e. $\ek=(1,\ldots,1,0,\ldots,0)^t$.
Clearly $\ek\not\in\LL\,.$ Hence
$\ZZ\cdot\ek\oplus\LL\subset\ZZ^N$ and
$(\ZZ\cdot\ek\oplus\LL)\otimes\CC^{\,\ast}\:\subset\:\CC^{\,\ast N}$.
Then
\begin{equation}\label{eq:pxtor}
\PP\EE_{\,\gT}\::=\:
\modquot{\CC^N_{\gT\circ}}{(\ZZ\cdot\ek\oplus\LL)\otimes\CC^{\,\ast}}\,.
\end{equation}
with the morphism $\PP\EE_{\,\gT}\rightarrow\PP_\gT$ induced from projection
onto the last $N-\kappa$ coordinates.

There are two kinds of codim $1$ simplices in the triangulation $\gT\,:$ those
which do contain $\core\gT$ and those which do not. Those which do not contain
$\core\gT$ are precisely the ones of the form $I\setminus\{i\}$ with
$I\in\gT^\na$ and $i\in\core\gT\,.$ Notice the relation with (\ref{eq:cnpt}).
The codim $1$ simplices which do not contain $\core\gT$ constitute a
triangulation of the boundary of $\Delta\,.$
Let as in (\ref{eq:a0})
\[
\ga_0\::=\:\ga_1+\ldots+\ga_\kappa\,.
\]
Then $\ZZ\cdot\ek\oplus\LL\subset\ZZ^N$ is precisely the lattice of
linear relations between the vectors
$\ga_1-\frac{1}{\kappa}\ga_0,\:\ga_2-\frac{1}{\kappa}\ga_0,\ldots,
\ga_N-\frac{1}{\kappa}\ga_0\,.$ Thus we see:

\

\begin{proposition}
$\PP\EE_{\,\gT}$ is the $(\na-1)$-dimensional smooth projective toric variety
associated with the lattice
$\ZZ (\ga_1-\frac{1}{\kappa}\ga_0)+\ldots+\ZZ(\ga_N-\frac{1}{\kappa}\ga_0)$ and
the fan consisting the cones with apex $0$ over the simplices of the
triangulation of the boundary of $-\frac{1}{\kappa}\ga_0+\Delta$ induced by
$\gT\,.$
\qed
\end{proposition}


\section{Calabi-Yau complete intersections in toric varieties}
\label{cicy}

According to proposition \ref{corecs} conditions (\ref{eq:core1}),
(\ref{eq:core2}), (\ref{eq:vol1}) imply that $\Lambda^\vee$
is a completely split reflexive Gorenstein cone. In \cite{babo}
Batyrev and Borisov relate this splitting property to complete intersections
in toric varieties. Formulated in our present context this relation
is as follows.

A (global) section of $\EE_{\,\gT}^\vee\rightarrow\PP_\gT$ is given by
polynomials $P_i(x_{\kappa+1},\ldots,x_N)$  $(i=1,\ldots,\kappa)$
which satisfy the homogeneity condition
\begin{equation} \label{eq:homogeneous}
P_i(t^{\ell_{\kappa+1}}\cdot x_{\kappa+1},\ldots,
t^{\ell_N}\cdot x_N)\:=\:t^{-\ell_i}\cdot P_i(x_{\kappa+1},\ldots,x_N)
\end{equation}
for every $ t\in\CC^{\,\ast} $ and $\ell=(\ell_1,\ldots,\ell_N)^t\in\LL\,.$
The vector bundle is a direct sum of line bundles and the polynomial $P_i$
gives a section of the $i$-th line bundle.

The polynomial $P_i$ is a linear combination of monomials
$x_{\kappa+1}^{m_{\kappa+1}}\cdot\ldots\cdot x_N^{m_N}$ such that
\[
\ell_{\kappa+1}m_{\kappa+1}+\ldots+\ell_N m_N\,=\,-\ell_i
\hspace{4mm}\textrm{ for all }\ell=(\ell_1,\ldots,\ell_N)\in\LL\,.
\]
These monomials correspond bijectively to the elements
$(m_1,\ldots,m_N)$ in the row space of matrix $\sA$ which satisfy
$m_i=1$, $m_j=0$ if $1\leq j\leq\kappa, j\neq i$ and $m_j\geq 0$ if
$j>\kappa\,.$ Equivalently, these monomials correspond bijectively to the
elements $\ew\in\ZZ^{\na\vee}$ which satisfy
\begin{equation}\label{eq:splitting2}
\ew\cdot\ga_j\:=\:
\left\{\begin{array}{ll}\geq  0 & \textrm{ for } j=1,\ldots,N\\
 1 & \textrm{ if } j=i\\
 0 & \textrm{ if } 1\leq j\leq\kappa\,,\;j\neq i
\end{array}\right.
\end{equation}
So the monomials in the polynomial $P_i$ correspond bijectively to the integral
lattice points in the polytope $\square_i$; see (\ref{eq:boxi}).

The zero locus of the section of $\EE_{\,\gT}^\vee\rightarrow\PP_\gT$
corresponding to the polynomials $P_i(x_{\kappa+1},\ldots,x_N)$
$(i=1,\ldots,\kappa)$ is clearly the complete intersection in $\PP_\gT$
with (homogeneous) equations
\begin{equation}\label{eq:cyci}
P_i(x_{\kappa+1},\ldots,x_N)=0\hspace{5mm}  (i=1,\ldots,\kappa)
\end{equation}
If the coefficients of these polynomials satisfy a
$\Lambda^\vee$-regularity condition, then this complete intersection
is a Calabi-Yau variety $Y$ of dimension $\na-2\kappa\,.$

The ring $\cR_{\sA,\gT}$ is isomorphic to the cohomology ring
of the toric variety $\PP_\gT\,.$ The elements $-c_1,\ldots,-c_{\kappa}$
are the Chern classes of the hypersurfaces associated with the
polynomials $P_1,\ldots,P_{\kappa}$.
With as before $c_{\core}=c_1\cdot\ldots\cdot c_{\kappa}$,
the ring $\modquot{\cR_{\sA,\gT}}{\Ann c_{\core}}$ is isomorphic
to the image of $H^\ast(\PP_\gT,\ZZ)$ in $H^\ast(Y,\ZZ)$.

\begin{center}\section*{Conclusions}\end{center}
Consider the map
$\ev:\,\CC^{N\vee}\rightarrow\CC^{N\vee}\:,$
$\;\ev(z_1,\ldots,z_N)\::=\:
(\ee^{2\pi i\,z_1},\ldots,\ee^{2\pi i\,z_N})\,.$
According to \cite{gkz2} p.304 cor.1.7  there is a vector $b\in\cC_\gT$
such that
\begin{equation}\label{eq:discno}
E_{\sA} (\ev(\ez))\neq 0\hspace{4mm}\textrm{for all}
\hspace{3mm}
\ez\in\CC^{N\vee} \hspace{3mm}\textrm{such that}\hspace{3mm}
 p\,(\Im\,\ez)\in b+\cC_\gT\,;
\end{equation}
here $p:\RR^{N\vee}\rightarrow\LL_\RR^\vee$ denotes
the surjection dual to the inclusion $\LL\hookrightarrow\ZZ^\na$.
This shows how one can replace the domain of definition
$\cV_\gT$ of the functions $\Psi_{\gT,\beta}$  (cf. (\ref{eq:defdom}))
by a slightly smaller domain $\cV^\prime_\gT$ such that
on  $\ev(\cV^\prime_\gT)$ the function $E_{\sA} $ is nowhere zero.
The $\cD$-module $ H^{\na}(\TTT\rel\tsZ_{\es-1})$
is therefore defined on $\ev(\cV^\prime_\gT)$;
cf. theorem \ref{hygeoDmod}.
Its pullback to $\cV^\prime_\gT$ is the
$\cD_\gT$-module $H^{\na}(\TTT\rel\tsZ_{\es-1})\otimes \cO_\gT\,,$
where $\cO_\gT$ denotes the ring of holomorphic functions on
$\cV^\prime_\gT$ and $\cD_\gT$ denotes the corresponding ring
of differential operators.

The functions $\Psi_{\gT,\beta}$
are also defined on the domain $\cV^\prime_\gT$ and
$$
\Psi_{\gT,\beta}\:\in\: \cR_{\sA,\gT}\otimes \cO_\gT\,.
$$
$\cR_{\sA,\gT}\otimes \cO_\gT$ is a $\cD_\gT$-module with $\cR_{\sA,\gT}$
as its group of horizontal sections.

The following theorem summarizes the results of this paper:

\

\begin{theorem}\label{mainthm}
Let $\{\ga_1,\ldots,\ga_N\}$ be a finite subset of $\ZZ^\na$ which
satisfies condition \ref{gkzcond}. Let
$\Lambda:=\RR_{\geq 0}\ga_1+\ldots+\RR_{\geq 0}\ga_N$ be the
associated Gorenstein cone and $\Delta:=\conv \{\ga_1,\ldots,\ga_N\}\,.$
\begin{enumerate}
\item
If there exists a unimodular regular triangulation of $\Delta$,
then condition \ref{aisall} is satisfied and
$\ga_1,\ldots,\ga_N$ generate $\ZZ^\na$, i.e.
\begin{equation}\label{eq:unimodex}
\Delta\cap\ZZ^\na=\{\ga_1,\ldots,\ga_N\}\hspace{4mm}\textrm{and}
\hspace{4mm}\MM=\ZZ^\na\,.
\end{equation}
\item
 For every unimodular regular triangulation $\gT$ there is an isomorphism
of $\cD_\gT$-modules on $\cV^\prime_\gT$ :
\begin{equation}\label{eq:isodmod}
 H^{\na}(\TTT\rel\tsZ_{\es-1})\otimes \cO_\gT\;\simeq\;
\cR_{\sA,\gT}\otimes \cO_\gT
\end{equation}
through which $\xi_0$ corresponds with $\Psi_{\gT,0}$.
More generally $\xi_\mu$ corresponds with $\Psi_{\gT,-\mu}$ if
$\mu\in\Lambda\cap\ZZ^\na\,.$
\item
In particular if $\Lambda$ is a reflexive Gorenstein cone
of index $\kappa$
and $\gT$ is a unimodular regular triangulation,
then $\cW_\na H^\na(\TTT\rel\tsZ_{\es-1})\otimes \cO_\gT$
is generated as a $\cD_\gT$-module by $\xi_{\ga_0}$ and
corresponds via (\ref{eq:isodmod}) with the sub-$\cD_\gT$-module
of $\cR_{\sA,\gT}\otimes \cO_\gT$ generated by $\Psi_{\gT,-\ga_0}\,.$

Moreover  $\xi_{\ga_0}$ has weight $\na$ and
Hodge type $(\na-\kappa,\kappa)$.
\item
If $\Lambda$ is a reflexive Gorenstein cone and $\gT$ is a unimodular regular
triangulation with non-empty core, then (\ref{eq:isodmod})
induces an isomorphism
\begin{equation}\label{eq:isodmodw}
\begin{array}{rcl}
\cW_\na  H^{\na}(\TTT\rel\tsZ_{\es-1})\otimes \cO_\gT &\simeq&
c_{\core}\cR_{\sA,\gT}\otimes \cO_\gT \\[.6em]
&\simeq& \modquot{\cR_{\sA,\gT}}{\Ann c_{\core}}\:\otimes \cO_\gT
\end{array}
\end{equation}
\item
Now assume $\gT$ satisfies  conditions
(\ref{eq:core1}), (\ref{eq:core2}), (\ref{eq:vol1}), i.e.
$\gT$ is a unimodular regular triangulation whose
core is not empty and is not contained in the boundary of $\Delta\,.$
Then
\begin{enumerate}
\item
$\Lambda$ is a reflexive Gorenstein cone.
\item
$\cR_{\sA,\gT}$ is isomorphic to the cohomology ring
$H^\ast(\PP_\gT,\ZZ)$ of the $(\na-\kappa)$-dimensional smooth
projective toric variety $\PP_\gT$:
\begin{equation}\label{eq:rattor}
\cR_{\sA,\gT}^{(m)}\simeq H^{2m}(\PP_\gT,\ZZ)\;,\hspace{5mm}
m=0,1,\ldots,\na-\kappa
\end{equation}
and in particular for $m=1$:
$\LL_\ZZ^\vee\simeq\textrm{Pic}(\PP_\gT)\,.$
\item
$c_{\core}=c_\kappa(\EE_{\,\gT})\,,$ the top Chern class
of the vectorbundle $\EE_{\,\gT}\,.$
\item
The zero locus of a general section of the dual vector bundle
$\EE_{\,\gT}^\vee$ is an $\na-2\kappa$-dimensional
Calabi-Yau complete intersection in $\PP_\gT$.
\item
\begin{eqnarray}
 H^{\na}(\TTT\rel\tsZ_{\es-1})\otimes \cO_\gT&\simeq&
H^\ast(\PP_\gT,\ZZ)\otimes \cO_\gT
\\
\cW_\na  H^{\na}(\TTT\rel\tsZ_{\es-1})\otimes \cO_\gT &\simeq&
c_\kappa(\EE_{\,\gT})\,H^\ast(\PP_\gT,\ZZ)\otimes \cO_\gT
\end{eqnarray}
\item
The monodromy representation is isomorphic to the representation
of $\textrm{Pic}(\PP_\gT)$ on $H^\ast(\PP_\gT,\ZZ)$
( resp. on $c_\kappa(\EE_{\,\gT})\,H^\ast(\PP_\gT,\ZZ)$ ) in which
the Chern class $c_1(\cL)$ of a line bundle $\cL$ acts as multiplication by
$\exp(c_1(\cL))\,.$
\end{enumerate}
\end{enumerate}
\end{theorem}
\textbf{Proof:}
\textbf{(i)}: corollary \ref{allvert}.
\textbf{(ii)}: theorems \ref{isosol} and \ref{hygeoDmod},
formulas (\ref{eq:recur}) and (\ref{eq:dximu2}).
\textbf{(iii)}: formulas (\ref{eq:alint}), (\ref{eq:a0z}),
(\ref{eq:weightfilt}) and theorem \ref{TZ}.
\textbf{(iv)}: corollary \ref{corid} and theorem \ref{corsol}.
\textbf{(v}\textit{a}\textbf{)}: proposition \ref{corecs}.
\textbf{(v}\textit{b}\textbf{)}: theorem \ref{protor} and
corollary \ref{allvert}.
\textbf{(v}\textit{c}\textbf{)}: section \ref{cicy}.
\textbf{(v}\textit{d}\textbf{)}: section \ref{cicy}.
\textbf{(v}\textit{e}\textbf{)}:
(\ref{eq:isodmod}) and (\ref{eq:isodmodw}).
\textbf{(v}\textit{f}\textbf{)}: formula (\ref{eq:monodromy}).
\qed

\

All cases which have on the A-side of mirror symmetry a smooth
complete intersection Calabi-Yau variety in a smooth projective
toric variety, are covered by this theorem. Indeed,
a smooth projective toric variety $\PP$ of dimension $d$ can be
constructed from a complete simplicial fan in which every maximal
cone is generated by a basis of the lattice $\ZZ^d$. Let
$\eu_1,\ldots,\eu_p\in\ZZ^d$ be the generators of the $1$-dimensional
cones in the fan and let
$$
\overline{\LL}:=\{(m_1,\ldots,m_p)\in\ZZ^p\:|\:m_1\eu_1+\ldots
m_p\eu_p=0\:\}
$$
The toric variety $\PP$ can also be obtained as the quotient of a
certain open part of $\CC^p$ by the action of the subtorus
$\overline{\LL}\otimes\CC^\ast$ of $(\CC^\ast)^p$. The Calabi-Yau
complete intersection $Y$ of codimension $\kappa$ in $\PP$ is the
common zero locus of polynomials $P_1,\ldots,P_\kappa$ which are
homogeneous for the action of $\overline{\LL}\otimes\CC^\ast$.
The homogeneity of $P_i$ is given by a character of this torus,
i.e. by a linear map $\chi_i:\overline{\LL}\rightarrow\ZZ\,.$
Now set $N=p+\kappa$ and $n=d+\kappa\,.$ Let
$$
\LL:=\{(-\chi_1(\sfm),\ldots,-\chi_\kappa(\sfm),m_1,\ldots,m_p)
\in\ZZ^N\:|\:\sfm=(m_1,\ldots,m_p)\in\overline{\LL}\:\}\,.
$$
Then $\LL$ has rank $N-n$. The Calabi-Yau condition for $Y$ implies
$\ell_1+\ldots+\ell_N=0$ for every
$\ell=(\ell_1,\ldots,\ell_N)\in\LL$.

Let $\sB$ be an $(N-n)\times N$-matrix
with entries in $\ZZ$ such that the columns of $\sB^t$ constitute
a basis for $\LL\,.$
Let $\sA$ be an $n\times N$-matrix of rank $n$
with entries in $\ZZ$ such that $\sA\cdot\sB^t=0\,.$
Then the columns $\ga_1,\ldots,\ga_N$ of $\sA$ satisfy condition
\ref{gkzcond}. One obtains a regular triangulation of
$\Delta:=\conv\{\ga_1,\ldots,\ga_N\}$ which satisfies the three conditions
(\ref{eq:core1}), (\ref{eq:core2}), (\ref{eq:vol1}),
by taking as its maximal simplices all
$\conv\{\ga_1,\ldots,\ga_\kappa,\ga_{\kappa+i_1},\ldots,\ga_{\kappa+i_d}\}$
for which $\eu_{i_1},\ldots,\eu_{i_d}$ span a maximal cone in the fan
defining $\PP$.

\subsection*{Acknowledgments}
I want to thank the Taniguchi Foundation and the organizers of the
Taniguchi Workshop ``Special Differential Equations'' 1991 for inviting
me to this interesting workshop where I first learnt about GKZ hypergeometric
functions and regular triangulations. Special thanks are for Bruce Hunt, who
kindly pointed out that the example in my talk \cite{sti} at the workshop
reminded him
very much of the examples
in Batyrev's talk a few weeks earlier in Kaiserslautern.

I want to thank the Japan Society for
the Promotion of Science for a JSPS Invitation Fellowship
in November-December 1996 and Kobe University
for support for a visit in July 1997.
The  stimulating atmosphere I experienced during these two visits to Kobe
was very important for finishing this work. I also want to thank
the organizers of various workshops - in particular
the Taniguchi Symposium ``Integrable Systems and Algebraic Geometry''
1997 - to which I was invited during these visits.
I am very grateful to Masa-Hiko Saito for arranging as my host both visits
and for the pleasant co-operation.


\end{document}